\documentclass[aos,seceqn,nameyear,dvips]{arximspdf}
\usepackage{graphics}
\doi{10.1214/09-AOS782}
\volume{38}
\issue{4}
\pubyear{2010}
\firstpage{2314}
\lastpage{2350}

\makeatletter

\newcommand{\mathds}{\mathbb}

\newtheorem{Theorem}{Theorem}[section]
\newtheorem{LEM}{Lemma}[section]
\newproclaim{EX}{Example}[section]
\newproclaim{RE}{Remark}

\makeatother

\begin{document}
\begin{frontmatter}

\title{Order thresholding\protect\thanksref{T1}}
\runtitle{Order thresholding}

\thankstext{T1}{Supported in part by NSF Grants
SES-0318200 and DMS-08-05598.}

\begin{aug}
\author[A]{\fnms{Min Hee} \snm{Kim}\ead[label=e1]{mzk132@psu.edu}} and
\author[A]{\fnms{Michael G.} \snm{Akritas}\corref{}\ead[label=e2]{mga@stat.psu.edu}}
\runauthor{M. H. Kim and M. G. Akritas}
\affiliation{Pennsylvania State University}
\address[A]{Department of Statistics\\
Pennsylvania State University\\
University Park, Pennsylvania 16802\\
USA\\
\printead{e1}\\
\phantom{E-mail: }\printead*{e2}}
\end{aug}

\received{\smonth{4} \syear{2009}}
\revised{\smonth{12} \syear{2009}}

%
\begin{abstract}
A new thresholding method, based on \textit{L}-statistics and called
\textit{order thresholding}, is proposed as a technique for improving
the power
when testing against high-dimensional alternatives. The new method
allows great flexibility in the choice of the threshold parameter. This
results in improved power over the soft and hard thresholding methods.
Moreover, order thresholding is not restricted
to the normal distribution. 
An extension of the basic order threshold statistic to high-dimensional
ANOVA is presented. The performance of the basic order threshold
statistic and its extension is evaluated with extensive simulations.
\end{abstract}

%

\end{frontmatter}

\section{Introduction}\label{sec1}

It is well known that, when testing against a high-dimen\-sional
alternative, omnibus tests designed to detect any departure from the
null hypothesis have low power. Neyman's (\citeyear{N37}) truncation idea,
though motivated by a different type of problem, served as the
spring board for the development of modern related approaches. Soft
and hard thresholding were introduced in the context of
nonparametric function estimation using wavelets by \citet{DJ94}. \citet{JS04} elaborate on a
number of additional applications of thresholding including image
processing, model selection, and data mining. \citet{B04}
considered applications to the one-way ANOVA design. \citet{S96}, \citet{F96} and \citet{FL98} consider applications of
thresholding methods to testing problems. \citet{F96} found that hard
thresholding outperforms both soft thresholding and adaptive
Neyman's truncation.

This paper proposes a new thresholding method based on $L$-statistics,
which is termed \textit{order
thresholding}. Order thresholding allows great flexibility in the
choice of the threshold parameter, can be
used for distributions other than the normal, and extends naturally to
factorial design settings.

In the simple context where the $X_i$ are independent $N(\theta_i,1)$,
$i=1,\ldots,n$, and we wish to test
$H_0\dvtx\theta_1=\cdots=\theta_n=0$ vs. $H_a\dvtx\theta_i\neq
0$, for
some $i$, the hard thresholding and order
thresholding test statistics are, respectively,
%
%
\begin{equation}\label{rel.t.s}
T_H(\delta_n)=\sum_{i=1}^nY_iI\{Y_i>\delta_n\}\quad \mbox{and}\quad
T_L(k_n)=\sum_{i=1}^nc_{in}Y_{i,n},
\end{equation}
where $Y_i=X_i^2$, $Y_{1,n}<\cdots<Y_{n,n}$ are the ordered $Y_i$'s,
$c_{in}=I(i> n-k_n)$, and $\delta_n$, $k_n$
are the corresponding threshold parameters. Thus, $T_L(k_n)$ is an
$L$-statistic based on the largest $k_n$
squared observations. Conceptually, the connection between hard
thresholding and order thresholding is similar
to that between type I and type II censoring. The main difference being
that the threshold parameters in type I
and type II censoring (the cut-off point and the proportion of
observations included, resp.) remain
fixed, while in the present case they change with the sample size. As
we will see, this distinction implies
very different asymptotic behavior.

The idea behind both statistics in (\ref{rel.t.s}) is similar to that
of Neyman's truncation. Namely, when the
``signal'' is known to be concentrated in a few locations, the
accumulation of stochastic errors has a negative
impact on the performance of the procedure based on the chi-square statistic
%
%
\begin{equation}\label{rel.naive.stat}
T_L(n)=\sum_{i=1}^nY_i.
\end{equation}
Since the signal locations are not known, the statistics in (\ref
{rel.t.s}) attempt to minimize the
accumulation of noise by focusing on the observations with the largest
absolute values. The asymptotic theory
for hard thresholding [\citet{F96}] requires several restrictive
conditions that prevent its general
applicability. For example, the centering and scaling of $T_H(\delta
_n)$ in (\ref{rel.t.s}) are specific to the
normality assumption and to the choice of $\delta_n$. Moreover,
$\delta_n$ is required to tend to infinity at a
rate that is specific to the normality assumption. For example, $\delta
_n$ tending to infinity is clearly not
appropriate if the $X_i$ have bounded support. (Below, we discuss an
application to multiple testing where the
$X_i$ are uniformly distributed.) Intuitively, if the signal is present
in more locations, it is advantageous
to lower the value of the hard threshold parameter. The advantage of
allowing different values of the threshold
parameter is amply illustrated in \citet{JS04}.
However, the asymptotic theory of
$T_H(\delta_n)$ requires the threshold parameter to tend to infinity
at specific rates. In particular, it must
be of the form $\delta_n = 2\log(na_n)$, where $a_n=c(\log n)^{-d}$,
for $c>0$ and $d>0.5$. Thus, if we let
$k_H(\delta_n)$ denote the random number of observations considered in
$T_H(\delta_n)$, the asymptotic theory
of $T_H(\delta_n)$ requires $E [k_H(\delta_n) ]$ to
converge to infinity at the rate of
\[
\frac{(\log n)^d}{\sqrt{\log n +d\log(c^{1/d}(\log n)^{-1})}}
\]
or, roughly, $(\log n)^{d-0.5}$. In contrast, the asymptotic theory of
$T_L(k_n)$ allows the threshold parameter
$k_n$ to tend to infinity at any rate.

While the asymptotic theory of $T_H(\delta_n)$ allows some flexibility
in the choice of~$\delta_n$, the
convergence of the distribution of $T_H(\delta_n)$ to its limiting
distribution is very slow unless $c=1$ and
$d=2$ [\citet{F96}]. The following tables show that small departures in
the recommended value of $d$, while
keeping $c=1$, have significant effect on the level of the test. The
results are based on 30,000 simulation
runs.

%
%
\begin{table}
\caption{Type \textup{I} errors of $T_H(\delta_n)$ for different values of the
hard threshold parameter}\label{paper1-t1}
\begin{tabular*}{\tablewidth}{@{\extracolsep{\fill}}lcccccc@{}}
\hline
& $\bolds{\delta_n-2.0}$ & $\bolds{\delta_n-1.6}$ & $\bolds{\delta_n-1.2}$
& $\bolds{\delta_n-0.8}$ & $\bolds{\delta_n-0.4}$ & $\bolds{\delta_n}$\\
\hline
$n=50$& 0.0003 & 0.0099 & 0.0231 & 0.0341 & 0.0431 & 0.0493\\
$n=100$& 0.0101 & 0.0229 & 0.0324 & 0.0390 & 0.0461 & 0.0504\\
$n=200$& 0.0231 & 0.0316 & 0.0382 & 0.0439 & 0.0484 & 0.0507\\
$n=500$& 0.0327 & 0.0388 & 0.0422 & 0.0465 & 0.0502 & 0.0535\\
\hline
\end{tabular*}
\end{table}

%
%
\begin{table}[b]
\caption{Type \textup{I} errors of $T_H(\delta_n)$ for different values of the
hard threshold parameter}\label{paper1-t2}
\begin{tabular*}{\tablewidth}{@{\extracolsep{\fill}}lccccc@{}}
\hline
& $\bolds{\delta_n+0.4}$ & $\bolds{\delta_n+0.8}$ & $\bolds{\delta_n+1.2}$
& $\bolds{\delta_n+1.6}$ & $\bolds{\delta_n+2.0}$\\
\hline
$n=50$& 0.0543 & 0.0588 & 0.0614 & 0.0654 & 0.0663\\
$n=100$& 0.0552 & 0.0559 & 0.0597 & 0.0616 & 0.0631\\
$n=200$& 0.0539 & 0.0562 & 0.0590 & 0.0601 & 0.0627\\
$n=500$& 0.0540 & 0.0563 & 0.0583 & 0.0604 & 0.0623\\
\hline
\end{tabular*}
\end{table}

To fully appreciate the results reported in Table \ref{paper1-t1}, we
mention that for $n=500$ the recommended
$\delta_{500}$ value is 5.1216, while the value $\delta_{500}-2$
corresponds to $c=1$ and $d=2.5474$. We see
that even with this small departure from the recommended value, the
achieved alpha level is 0.0327 even with
$n=500$. To contrast these results with those of Table \ref{paper1-t3} for
$T_L(k_n)$, note that in Table \ref{paper1-t3} with $n=500$,
$k_n$ ranges from 2 to 500, while in Table \ref{paper1-t1} with
$n=500$, the $E [k_H(\delta_n+h) ]$
ranges from 38.63 for $h=-2.0$ to 11.81 for $h=0$. Thus, the
deterioration of the achieved alpha levels occurs
as $E [k_H(\delta_n) ]$ increases over a relatively small
range (in each case, the variance of
$k_H(\delta_n)$ is slightly smaller than its expected value). In Table
\ref{paper1-t2} with $n=500$, the
$E [k_H(\delta_n+h) ]$ ranges from 9.39 for $h=0.4$ to 3.80
for $h=2.0$, and for this range of values
the type I error rate does not change much. In both tables with
$n=500$, the variance of the binomial random
variable $k_H(\delta_n+h)$ is slightly smaller than its expected value
because $P(Y_i\leq\delta_n+h)>0.92$ for
$-2.0\leq h\leq2.0$. Finally, following a remark by the AE, we note
that the slightly liberal $\alpha$ levels
of the $T_L(k_n)$ statistic can be corrected by the use of a multiple
of a $\chi^2$ distribution to approximate
its finite sample distribution. Thus, using the approximation
$T_L(k_n)\stackrel{\cdot}{\sim} b\chi^2_\nu$, where
$b$ and $\nu$ are chosen to match the mean and variance of $T_L(k_n)$,
results in the type I
error rates shown in Table \ref{paper1-t3-app}.

%
%
\begin{table}
\caption{Type \textup{I} errors of $T_L(k_n)$ for different values of the order
threshold parameter}\label{paper1-t3}
\begin{tabular*}{\tablewidth}{@{\extracolsep{\fill}}lcccccccc@{}}
\hline
& $\bolds{[\log^{1/2} n]}$ & $\bolds{[\log n]}$ & $\bolds{[\log^{3/2} n]}$
& $\bolds{[n^{1/2}]}$ & $\bolds{[n^{2/3}]}$ &
$\bolds{[n^{3/4}]}$ & $\bolds{[n^{7/8}]}$ & $\bolds{n}$\\
\hline
$n=50$& 0.0696& 0.0685& 0.0646& 0.0646& 0.0635& 0.0630& 0.0623& 0.0626\\
$n=100$& 0.0669& 0.0640& 0.0606& 0.0600& 0.0591& 0.0577& 0.0585&
0.0589\\
$n=200$& 0.0667& 0.0620& 0.0603& 0.0589& 0.0582& 0.0583& 0.0555&
0.0560\\
$n=500$& 0.0665& 0.0631& 0.0577& 0.0559& 0.0536& 0.0536& 0.0535&
0.0547\\
\hline
\end{tabular*}
\end{table}

%
%
\begin{table}[b]
\caption{Type \textup{I} error rates using the approximation $T_L(k_n)\stackrel
{\cdot}{\sim} b\chi^2_\nu$}\label{paper1-t3-app}
\begin{tabular*}{\tablewidth}{@{\extracolsep{\fill}}lcccccccc@{}}
\hline
& $\bolds{[\log^{1/2} n]}$ & $\bolds{[\log n]}$ & $\bolds{[\log^{3/2} n]}$ & $\bolds{[n^{1/2}]}$
& $\bolds{[n^{2/3}]}$ & $\bolds{[n^{3/4}]}$ & $\bolds{[n^{7/8}]}$ & $\bolds{n}$\\
\hline
$n=50$& 0.0565& 0.0545& 0.0531& 0.0531& 0.0532& 0.0534& 0.0540& 0.0546\\
$n=100$& 0.0566& 0.0540& 0.0520 & 0.0522& 0.0519& 0.0507& 0.0518&
0.0520\\
$n=200$& 0.0555& 0.0547& 0.0536 & 0.0531& 0.0526& 0.0523& 0.0530&
0.0516\\
$n=500$& 0.0589& 0.0556& 0.0552 & 0.0530& 0.0520& 0.0521& 0.0508&
0.0505\\
\hline
\end{tabular*}
\end{table}

\begin{figure}

\includegraphics{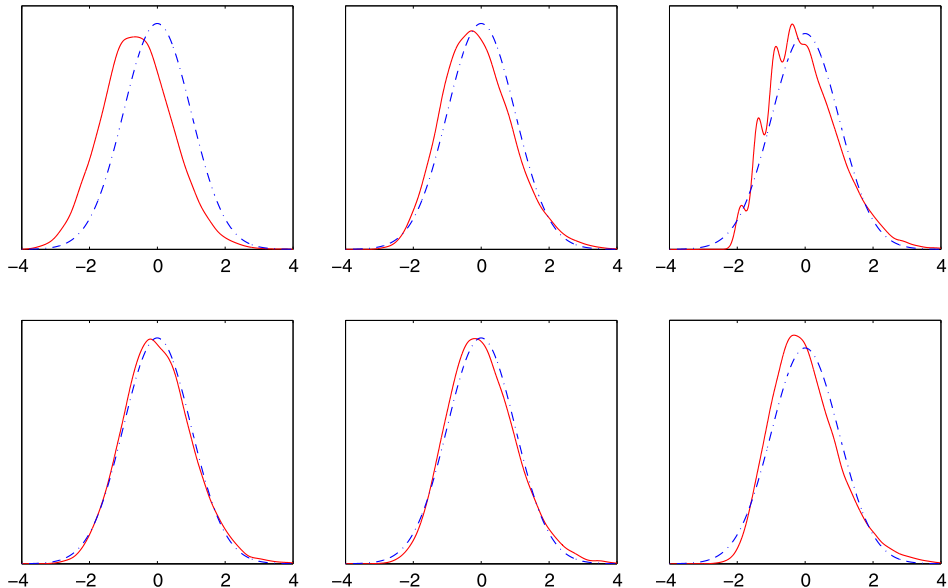}

\caption{Top panel: estimated densities of $T_H(\delta_{200})$ for
$\delta_{200}= 1.842$, 3.927 and 5.672.
Bottom panel: estimated densities of $T_L(k_{200})$ for $k_{200} = 35$,
10 and 3.}
\label{fig.exponential_on_normal_pp}
\end{figure}

The greater flexibility in the choice of the threshold parameter that
the order threshold statistic offers does
not come at the expense of the rate with which it converges to its
asymptotic distribution. To emphasize this
aspect, Figure \ref{fig.exponential_on_normal_pp} presents the
estimated densities of the hard thresholding
(solid lines in the upper panel) and order thresholding test statistics
(solid lines in the lower panel), based
on 20,000 simulated values of each statistic using $n=200$. The
threshold parameters of the hard thresholding
and order thresholding test statistics have been chosen so that the
average number of observations included in
the two statistics are the same in each column. We see that the
estimated densities of the order threshold
statistic are closer to the standard normal density (dash-dot line)
than those of the hard threshold statistic.
In particular, the estimated densities in the upper panel show the
rapid deterioration of the quality of the
normal approximation to the distribution of $T_H(\delta_{200})$ as
$\delta_{200}$ shifts away from recommended
value of $\delta_{200} = 2\log(200\log^{-2} 200)=3.9271$.



%
%

The remaining sections of this paper are organized as follows. In
Section \ref{sec2},
we represent a special form of the order statistics using data from
an exponential distribution and briefly review the methodology of
\citet{CGJ67}. Section \ref{sec3} develops the order
threshold procedure for testing normal means in settings where the
number of parameters increases with the sample size, presents
simulation results comparing the hard thresholding, a power-enhanced
version of the \citet{S86}, and order thresholding test statistics,
and gives a recommendation for choosing the data-driven value of the
threshold parameter. Section \ref{sec4} extends the order thresholding test
procedure to the high-dimensional ANOVA setting [called HANOVA in
\citet{FL98}], presents simulation results comparing the power
of the classical $F$ and order threshold statistics, and gives a
recommendation for a data-driven choice of the order threshold
parameter. A discussion summarizing the developments is given in
Section \ref{sec5}. Finally, the condensed proofs are given in the Appendix.
For detailed proofs, see the archived supplemental
material in \citet{KA09}. This is part of the Ph.D. dissertation
of the first author.


\section{From order statistics to order thresholding: An overview}\label{sec2}

In the late 1960s when the asymptotic theory of linear combinations of
order statistics ($L$-statistics) was
developed [cf. \citet{B67}, \citet{CGJ67},
\citet{Shorack69}, \citet{Stigler69}] the main
emphasis was in the estimation of the location parameter. Therefore,
the conditions in these papers do not
yield automatically the asymptotic distribution of $L$-statistics that
assign positive weight to only the
largest order statistics. Such $L$-statistics were considered by
\citet{N82} in his study of the selection
differential for applications to outlier detection. Using results from
\citet{H78} and \citet{S73}, he
obtained the asymptotic distribution in the extreme and quantile cases,
respectively. Here, we will use the
conditions from the paper of \citet{CGJ67},
CGJ1967 from now on. Their approach is
based on a special representation of the order statistics from the
exponential distribution, which we now
review.

Let $V_1,\ldots,V_n$ be i.i.d. from the standard exponential
distribution, let $V_{1,n}<\cdots<V_{n,n}$ be the
corresponding order statistics, and consider the order threshold statistic
%
%
\begin{equation}\label{paper1-2.3}
T_{E,L}(k_n)=\sum_{i=1}^n c_{in}V_{i,n}=\sum_{i=n-k_n+1}^n V_{i,n}.
\end{equation}
The method of CGJ1967 for establishing the asymptotic distribution of
$T_{E,L}(k_n)$ rests on the following
well-known property [cf. \citet{DN03}, pages~17 and 18].
\begin{LEM} The vector of order statistics
$(V_{1,n},\ldots,V_{n,n})$ may be represented in distribution by
\[
(V_{1,n},\ldots,V_{n,n})\stackrel{d}{=} (Y_1,\ldots,Y_n),
\]
where
\[
Y_i=\frac{V_1}{n}+\frac{V_2}{n-1}+\cdots+\frac{V_i}{n-i+1} =\sum
_{j=1}^i\frac{V_j}{n-j+1}.
\]
Thus, with $T_{E,L}(k_n)$ given by (\ref{paper1-2.3}), it can be
represented in distribution as
%
%
\begin{equation} \label{paper1-2.5}
T_{E,L}(k_n)\stackrel{d}{=}
\sum_{j=1}^n\alpha_{E,jn}(k_n)V_j,
\end{equation}
where 
$\alpha_{E,jn}(k_n)=k_n/(n-j+1)$ for $j\leq n-k_n$ and $\alpha
_{E,jn}(k_n)=1$ for $j>n-k_n$.
\end{LEM}

Relation (\ref{paper1-2.5}) expresses $T_{E,L}(k_n)$ as a linear
combination of the independent random variables $V_1,\ldots,V_n$
which enables the use of standard asymptotic results for
establishing conditions for its asymptotic distribution. This is
given, without proof, in the following.
\begin{Theorem}\label{th.exp.case}
Let $k_n, n\ge1$, be any sequence of integers which satisfies
$k_n\to\infty$, as $n\to\infty$, and $k_n\le n$, and let
$T_{E,L}(k_n)$ be given in (\ref{paper1-2.3}). Then we have
%
%
\begin{equation}\label{paper1-2.6}\quad
T_{E,L}^*(k_n)=\frac{T_{E,L}(k_n)-\sum_{i=1}^n\alpha_{E,in}(k_n)}
{\sqrt{\sum_{i=1}^n\alpha_{E,in}(k_n)^2}}
\stackrel{d}{\rightarrow}N(0,1)\qquad \mbox{as } n\to\infty.
\end{equation}
\end{Theorem}

In the case where the observations $Y_i$, $i=1,\ldots,n$, come from a
distribution function $F$, the CGJ1967
approach for obtaining the asymptotic distribution of the order
threshold statistic
%
%
\begin{equation}\label{paper1-2.3.1}
T_{F,L}(k_n)=\sum_{i=1}^n c_{in}Y_{i,n}=\sum_{i=n-k_n+1}^n Y_{i,n},
\end{equation}
where $Y_{1,n}<\cdots<Y_{n,n}$ are the ordered $Y_i$'s, is based on the
expression $Y_{i,n}=\tilde
H_F(V_{i,n})$, where $\tilde H_F= F^{-1}\circ G$, and $G$ is the
standard exponential distribution function,
and the use of Taylor expansion to obtain:
\begin{LEM}[{[\citet{CGJ67}]}]
Let $T_{F,L}(k_n)$ be given by (\ref{paper1-2.3.1}). Then,
\[
n^{-1}T_{F,L}(k_n)\stackrel{d}{=}\mu_{F,n}(k_n)+Q_{F,n}(k_n)+R_{F,n}(k_n),
\]
where
\begin{eqnarray*}
\mu_{F,n}(k_n)&=&\frac{1}{n}\sum_{i=1}^n c_{in}\tilde{H}_F
(\tilde{\nu}_{in} ),
\\
Q_{F,n}(k_n)&=&\frac{1}{n}\sum_{i=1}^n\alpha_{F,in}(k_n)
(V_i-1 )
\end{eqnarray*}
and
\[
R_{F,n}(k_n)=\frac{1}{n}\sum_{i=1}^nc_{in} \bigl\{ \bigl(\tilde
{H}_F (V_{i,n} )
-\tilde{H}_F (\tilde{\nu}_{in} ) \bigr)-
(V_{i,n}-\tilde{\nu}_{in} )\tilde{H}_F'
(\tilde{\nu}_{in} ) \bigr\}
\]
with $\alpha_{F,in}(k_n)=\frac{1}{n-i+1}
\sum_{j=i}^nc_{jn}\tilde{H}_F' (\tilde{\nu}_{jn} )$ and
$\tilde{\nu}_{in}=\sum_{j=1}^i\frac{1}{n-j+1}$.
\end{LEM}

They then provide conditions under which $Q_{F,n}(k_n)$ is
asymptotically normally distributed and the
remainder term, $R_{F,n}(k_n)$, tends to zero in probability.

\section{Single sequence of $N(0,1)$ random variables}\label{sec3}

In this section, we will apply the approach of CGJ1967 to develop order
threshold test procedures for testing
the simple hypothesis
%
%
\begin{equation}\label{rel.nullH.sing.seq}
H_0\dvtx\theta_i=0\qquad \forall i \quad\mbox{versus}\quad H_a\dvtx\mbox{$H_0$ is false}
\end{equation}
based on a sequence of observations $X_i, i=1,\ldots,n$, where
$X_i\sim N(\theta_i,1)$. The asymptotic null
distribution of the order threshold statistic given by (\ref{rel.t.s})
is derived in the next subsection, while
simulation results comparing the power of the hard threshold statistic,
a power-enhanced version of the \citet{S86} statistic, and that of order threshold statistics are presented
in Section \ref{ss.sim.sing.seq}. The
simulation results suggest that choosing the order threshold parameter
equal to the number of the false null
hypotheses maximizes the power. Section \ref{choosingkn} presents a
recommendation for a data-driven choice
of the order threshold parameter using the idea of Storey (\citeyear{S02}, \citeyear{S03}).

\subsection{The asymptotic null distribution}\label{sec31}

Let $X_i, i=1,\ldots,n$, be standard normal random variables, and let
%
%
\begin{equation}\label{paper1-2.7}
T_L(k_n)=\sum_{i=1}^nc_{in}Y_{i,n}=\sum_{i=n-k_n+1}^n Y_{i,n},
\end{equation}
where $Y_i=X_i^2$, $Y_{1,n}<\cdots<Y_{n,n}$ are the ordered $Y_i$'s,
$c_{in}=I(i>n-k_n)$, and $k_n$ is the order
threshold parameter. The approach of CGJ1967 is based on the representation
\[
T_L(k_n)\stackrel{d}{=}\sum_{i=n-k_n+1}^n\tilde{H}
(V_{i,n} ),
\]
where $V_{i,n}$, $i=1,\ldots,n$, are the ordered observations from an
i.i.d. sequence of $\operatorname{Exp}(1)$ random
variables, and
\[
\tilde{H}(v)=F^{-1}\circ G(v)
\]
with $F(y)=\frac{1}{\sqrt{2\pi}} {\int
_{0}^{y}}t^{-1/2}e^{-t/2}\,dt, y>0$, and
$G(v)=1-e^{-v}$, $v\geq0$. Let
%
%
\begin{equation} \label{eq:1}
\mu_n(k_n)=\frac{1}{n}\sum_{i=1}^n c_{in}\tilde{H} (\tilde
{\nu}_{in} ),\qquad
\sigma_n^2(k_n)=\frac{1}{n}\sum_{i=1}^n\alpha_{in}^2(k_n),
\end{equation}
where
%
%
\begin{eqnarray}\label{alpha.tilednu}
\alpha_{in}(k_n)&=&\frac{1}{n-i+1} \sum_{j=i}^nc_{jn}\tilde{H}'
(\tilde{\nu}_{jn} ),\nonumber\\[-8pt]\\[-8pt]
\tilde{\nu}_{in}&=&\sum_{j=1}^i\frac{1}{n-j+1}.\nonumber
\end{eqnarray}
The term of $\alpha_{in}(k_n)$ can be re-expressed as $\alpha
_{in}(k_n)=\frac{1}{n-i+1}
\sum_{j=n-k_n+1}^n\tilde{H}' (\tilde{\nu}_{jn} )$ for
$i\leq n-k_n$ and
$\alpha_{in}(k_n)=\frac{1}{n-i+1} \sum_{j=i}^n\tilde{H}'
(\tilde{\nu}_{jn} )$ for $i>n-k_n$ with
$\tilde{H}' (\tilde{\nu}_{jn} )=\frac{e^{-\tilde{\nu
}_{jn}}}{f (F^{-1}
(1-e^{-\tilde{\nu}_{jn}} ) )}$ and the function $f$
is the derivative of $F$. With this notation
we have the following.
\begin{Theorem}\label{sing.seq.3.1}
Let $Y_i$, $i=1,\ldots,n$, be a sequence of i.i.d. random variables
having the central chi-squared distribution with 1 degree of
freedom. Let $k_n, n\ge1$, be any sequence of integers which
satisfies $k_n\to\infty$, as $n\to\infty$, and $k_n\le n$. Let
$\mu_n(k_n)$ and $\sigma_n^2(k_n)$ be as in (\ref{eq:1}) with
$c_{in}=I (i>n-k_n )$, and let $T_L(k_n)$ be given in
(\ref{paper1-2.7}). Then we have
%
%
\begin{equation}\label{paper1-2.14}
T_L^*(k_n)=\frac{T_L(k_n)-n\mu_n(k_n)}{\sqrt{n}\sigma
_n(k_n)}\stackrel{d}{\rightarrow} N(0,1)\qquad
\mbox{as } n\to\infty.
\end{equation}
\end{Theorem}

Note that the asymptotic mean of $T_L(k_n)$ is $n\mu_n(k_n)$ and the
asymptotic variance of $T_L(k_n)$ is
$n\sigma_n^2(k_n)$ as $k_n$ tends to infinity with $n$.

\subsection{Simulations}\label{ss.sim.sing.seq}

In this subsection, we compare the empirical power of the order
threshold statistic using several values of the
threshold parameter with those of the hard threshold and a
power-enhanced version of \citet{S86}
statistics. The original Simes multiple testing procedure rejects the
global hypothesis, $H_{0}^G$, that all
$H_{0}^{(i)}\dvtx\theta_i=0$, $i=1,\ldots,n$, are true if
\[
T_S=\min_{1\leq i\leq n} \bigl\{nP_{(i)}/i \bigr\}<\alpha,
\]
where $P_{(1)}<\cdots<P_{(n)}$ are the ordered $p$-values of the
individual hypotheses, and $\alpha$ is the
desired level of significance. A power-enhanced version of the original
Simes test procedure uses
$\alpha/(1-k_n^{\mathrm{opt}}/n)$ instead of $\alpha$, where $k_n^{\mathrm{opt}}$ is
the number of false null hypotheses.

The simulations reported here use samples of size $n=500$ generated
from the normal distribution with variance
1. The threshold parameter $k_{500}$ of the order threshold statistics
takes values of 15, 40, 70, 100, 200,
500, as well as a data-driven value, denoted by $\widehat
k_{500}^{\mathrm{opt}}$, whose description is given in
Section \ref{choosingkn}. The empirical power using the
approximation $T_L(\widehat
k_{500}^{\mathrm{opt}})\stackrel{\cdot}{\sim} b\chi^2_\nu$ is reported
together with that using the normal approximation to
$T_L(\widehat k_{500}^{\mathrm{opt}})$. The hard threshold statistic we consider
uses the recommended value of the
threshold parameter which is $\delta_{500}= 2\log(500\log^{-2}
500)=5.1216$. All results are based on 3000
simulation runs. Since the the global hypothesis $H_{0}^G$ is the same
for all three simulation settings, the
type I error rates reported in the last row of Table \ref{paper1-sing-seq-t4} pertain also to
Tables~\ref{paper1-sing-seq-t6} and \ref{table7}. Note that all achieved
%
%
%
\begin{table}
\tabcolsep=0pt
\caption{Power calculations in Example
\protect\ref{sing.seq.ex4}}\label{paper1-sing-seq-t4}
\begin{tabular*}{\tablewidth}{@{\extracolsep{\fill}}lccccccccccc@{}}
\hline
& $\bolds{k_{500}^{\mathrm{opt}}}$ & $\bolds{T_S}$ &$\bolds{T_H(5.122)}$
& $\bolds{T_L(\widehat k_{500}^{\mathrm{opt}})}$ & $\bolds{b\chi_{\nu}^2}$
& $\bolds{T_L(15)}$ & $\bolds{T_L(40)}$ & $\bolds{T_L(70)}$
& $\bolds{T_L(100)}$ & $\bolds{T_L(200)}$ & $\bolds{T_L(500)}$\\
\hline
$H_{1}$& 30& 0.843& 0.944& 0.977& 0.975& 0.976& 0.973& 0.968& 0.960&
0.938& 0.913\\
$H_{3}$& 28& 0.845& 0.942& 0.978& 0.975& 0.976& 0.975& 0.969& 0.961&
0.937& 0.910\\
$H_{5}$& 26& 0.840& 0.926& 0.972& 0.970& 0.971& 0.966& 0.956& 0.943&
0.911& 0.879\\
$H_{7}$& 24& 0.796& 0.893& 0.950& 0.948& 0.949& 0.942& 0.929& 0.915&
0.880& 0.851\\
$H_{8}$& 23& 0.777& 0.845& 0.933& 0.928& 0.932& 0.915& 0.891& 0.875&
0.818& 0.775\\
$H_{10}$& 21& 0.764& 0.817& 0.908& 0.900& 0.907& 0.891& 0.868& 0.841&
0.785& 0.744\\
$H_{11}$& 20& 0.766& 0.792& 0.905& 0.899& 0.906& 0.883& 0.853& 0.832&
0.764& 0.712\\
$H_{12}$& 19& 0.764& 0.783& 0.903& 0.897& 0.903& 0.873& 0.841& 0.812&
0.751& 0.709\\
$H_{13}$& 18& 0.750& 0.752& 0.881& 0.875& 0.880& 0.845& 0.804& 0.776&
0.709& 0.662\\
$H_{14}$& 17& 0.739& 0.734& 0.864& 0.858& 0.869& 0.836& 0.789& 0.760&
0.694& 0.649\\
$H_{15}$& 16& 0.559& 0.574& 0.724& 0.707& 0.723& 0.671& 0.633& 0.608&
0.541& 0.495\\
$H_{16}$& 15& 0.526& 0.564& 0.707& 0.693& 0.707& 0.660& 0.611& 0.574&
0.517& 0.484\\
$H_{17}$& 14& 0.532& 0.529& 0.675& 0.661& 0.677& 0.625& 0.574& 0.542&
0.467& 0.432\\
$H_{18}$& 13& 0.464& 0.435& 0.584& 0.568& 0.590& 0.534& 0.496& 0.458&
0.404& 0.373\\
$H_{19}$& 12& 0.483& 0.402& 0.570& 0.556& 0.574& 0.500& 0.459& 0.427&
0.374& 0.347\\
$H_{20}$& 11& 0.470& 0.380& 0.547& 0.533& 0.551& 0.475& 0.425& 0.395&
0.343& 0.308\\
$H_{21}$& 10& 0.467& 0.390& 0.555& 0.540& 0.559& 0.490& 0.433& 0.402&
0.341& 0.319\\
$H_{22}$& \phantom{0}9& 0.460& 0.364& 0.534& 0.515& 0.535& 0.454& 0.402& 0.368&
0.313& 0.281\\
$H_{23}$& \phantom{0}8& 0.460& 0.362& 0.517& 0.503& 0.522& 0.447& 0.389& 0.351&
0.301& 0.279\\
$H_{24}$& \phantom{0}7& 0.417& 0.290& 0.450& 0.434& 0.455& 0.375& 0.318& 0.288&
0.248& 0.230\\
$H_{0}^G$&\phantom{0}0&0.052 &0.050 &0.059 &0.057 & 0.057&
0.054&0.052&0.051&0.052& 0.055 \\
\hline
\end{tabular*}
\end{table}
%
%
%
\begin{table}
\tabcolsep=0pt
\caption{Power calculations in
Example \protect\ref{sing.seq.ex6}}\label{paper1-sing-seq-t6}
\begin{tabular*}{\tablewidth}{@{\extracolsep{\fill}}lccccccccccc@{}}
\hline
& $\bolds{k_{500}^{\mathrm{opt}}}$ & $\bolds{T_S}$ &$\bolds{T_H(5.122)}$
& $\bolds{T_L(\widehat k_{500}^{\mathrm{opt}})}$ & $\bolds{b\chi_{\nu}^2}$
& $\bolds{T_L(15)}$ & $\bolds{T_L(40)}$ & $\bolds{T_L(70)}$ &
$\bolds{T_L(100)}$ & $\bolds{T_L(200)}$ & $\bolds{T_L(500)}$\\
\hline
$H_{1}$& 30& 0.650& 0.574& 0.759& 0.745& 0.760& 0.699& 0.651& 0.608&
0.548& 0.513\\
$H_{2}$& 29& 0.680& 0.584& 0.755& 0.741& 0.761& 0.700& 0.649& 0.612&
0.544& 0.504\\
$H_{3}$& 28& 0.652& 0.565& 0.745& 0.729& 0.747& 0.684& 0.640& 0.602&
0.540& 0.498\\
$H_{6}$& 25& 0.666& 0.549& 0.728& 0.717& 0.732& 0.667& 0.625& 0.591&
0.521& 0.479\\
$H_{7}$& 24& 0.677& 0.562& 0.743& 0.729& 0.745& 0.686& 0.632& 0.591&
0.522& 0.482\\
$H_{8}$& 23& 0.666& 0.536& 0.716& 0.703& 0.724& 0.657& 0.612& 0.569&
0.508& 0.478\\
$H_{10}$& 21& 0.340& 0.350& 0.449& 0.434& 0.445& 0.418& 0.394& 0.367&
0.333& 0.317\\
$H_{12}$& 19& 0.351& 0.342& 0.444& 0.426& 0.443& 0.410& 0.383& 0.362&
0.341& 0.305\\
$H_{13}$& 18& 0.342& 0.330& 0.456& 0.442& 0.450& 0.416& 0.388& 0.367&
0.335& 0.316\\
$H_{14}$& 17& 0.350& 0.331& 0.448& 0.432& 0.451& 0.412& 0.377& 0.363&
0.325& 0.300\\
$H_{15}$& 16& 0.337& 0.334& 0.432& 0.416& 0.431& 0.402& 0.375& 0.356&
0.327& 0.307\\
$H_{16}$& 15& 0.330& 0.294& 0.406& 0.393& 0.403& 0.371& 0.338& 0.319&
0.293& 0.274\\
$H_{17}$& 14& 0.357& 0.282& 0.399& 0.387& 0.403& 0.352& 0.323& 0.305&
0.267& 0.252\\
$H_{18}$& 13& 0.325& 0.290& 0.393& 0.378& 0.390& 0.358& 0.329& 0.312&
0.276& 0.261\\
$H_{19}$& 12& 0.337& 0.296& 0.413& 0.396& 0.412& 0.368& 0.337& 0.314&
0.277& 0.255\\
$H_{20}$& 11& 0.343& 0.291& 0.399& 0.383& 0.399& 0.349& 0.314& 0.296&
0.270& 0.250\\
$H_{21}$& 10& 0.346& 0.290& 0.405& 0.391& 0.404& 0.356& 0.321& 0.306&
0.268& 0.248\\
$H_{22}$& \phantom{0}9& 0.224& 0.198& 0.264& 0.251& 0.262& 0.237& 0.220& 0.208&
0.195& 0.189\\
$H_{23}$& \phantom{0}8& 0.196& 0.190& 0.257& 0.242& 0.253& 0.228& 0.216& 0.197&
0.191& 0.182\\
$H_{24}$& \phantom{0}7& 0.207& 0.182& 0.256& 0.245& 0.253& 0.225& 0.212& 0.200&
0.186& 0.179\\
\hline
\end{tabular*}
\end{table}
%
%
\begin{table}
\tabcolsep=0pt
\caption{Power calculations in Example
\protect\ref{sing.seq.ex7}}\label{table7}
\begin{tabular*}{\tablewidth}{@{\extracolsep{\fill}}lccccccccccc@{}}
\hline
& $\bolds{k_{500}^{\mathrm{opt}}}$ & $\bolds{T_S}$ &$\bolds{T_H(5.122)}$
& $\bolds{T_L(\widehat k_{500}^{\mathrm{opt}})}$ & $\bolds{b\chi_{\nu}^2}$
& $\bolds{T_L(15)}$ & $\bolds{T_L(40)}$ & $\bolds{T_L(70)}$
& $\bolds{T_L(100)}$ & $\bolds{T_L(200)}$ & $\bolds{T_L(500)}$\\
\hline
$H_{1}$& 30& 0.674& 0.959& 0.973& 0.970& 0.969& 0.976& 0.982& 0.981&
0.970& 0.962\\
$H_{3}$& 28& 0.643& 0.954& 0.966& 0.960& 0.955& 0.974& 0.974& 0.973&
0.960& 0.947\\
$H_{4}$& 27& 0.617& 0.935& 0.957& 0.954& 0.945& 0.965& 0.963& 0.961&
0.950& 0.934\\
$H_{6}$& 25& 0.598& 0.900& 0.936& 0.931& 0.926& 0.947& 0.943& 0.941&
0.922& 0.903\\
$H_{8}$& 23& 0.566& 0.872& 0.912& 0.905& 0.902& 0.920& 0.917& 0.911&
0.891& 0.865\\
$H_{10}$& 21& 0.529& 0.831& 0.877& 0.869& 0.862& 0.889& 0.886& 0.875&
0.849& 0.817\\
$H_{12}$& 19& 0.509& 0.777& 0.837& 0.828& 0.821& 0.843& 0.833& 0.821&
0.786& 0.753\\
$H_{13}$& 18& 0.481& 0.740& 0.816& 0.803& 0.802& 0.813& 0.803& 0.785&
0.738& 0.703\\
$H_{14}$& 17& 0.472& 0.715& 0.785& 0.773& 0.772& 0.784& 0.773& 0.763&
0.710& 0.671\\
$H_{15}$& 16& 0.448& 0.674& 0.748& 0.732& 0.736& 0.749& 0.735& 0.715&
0.669& 0.633\\
$H_{16}$& 15& 0.418& 0.630& 0.715& 0.700& 0.702& 0.706& 0.686& 0.668&
0.624& 0.585\\
$H_{17}$& 14& 0.393& 0.569& 0.658& 0.645& 0.645& 0.646& 0.629& 0.610&
0.562& 0.523\\
$H_{18}$& 13& 0.368& 0.522& 0.629& 0.616& 0.623& 0.620& 0.597& 0.573&
0.523& 0.489\\
$H_{19}$& 12& 0.341& 0.498& 0.593& 0.577& 0.582& 0.582& 0.552& 0.525&
0.486& 0.451\\
$H_{20}$& 11& 0.328& 0.441& 0.539& 0.527& 0.539& 0.519& 0.491& 0.472&
0.436& 0.407\\
$H_{21}$& 10& 0.306& 0.390& 0.487& 0.470& 0.480& 0.464& 0.436& 0.421&
0.382& 0.353\\
$H_{22}$& \phantom{0}9& 0.285& 0.354& 0.439& 0.423& 0.438& 0.422& 0.393& 0.379&
0.344& 0.317\\
$H_{23}$& \phantom{0}8& 0.260& 0.298& 0.393& 0.374& 0.386& 0.367& 0.342& 0.318&
0.292& 0.276\\
$H_{24}$& \phantom{0}7& 0.245& 0.265& 0.349& 0.333& 0.346& 0.315& 0.296& 0.283&
0.255& 0.236\\
$H_{25}$& \phantom{0}6& 0.221& 0.224& 0.300& 0.286& 0.295& 0.272& 0.257& 0.245&
0.228& 0.213\\
\hline
\end{tabular*}
\end{table}
significance levels are below 0.06. The alternatives considered have 30
of the 500 mean values different from
zero. In particular, we consider the following sequence of alternatives
indexed by $r$:
\[
H_r\dvtx\theta_j=\eta_{j+r-1} \qquad\mbox{for } j=1,\ldots,500,
r=1,\ldots,30,
\]
where $\eta_j$, $j=1,2,\ldots,$ is a given sequence. The following
are examples with different values of
$\bolds\eta$.
\begin{EX}\label{sing.seq.ex4}
We generate the values of $\eta_{j}$, $j=1,\ldots,30$, from
$N(1.5,1)$. The rest values of $\eta_j$ are 0. The
values different from 0 are as follows:
\begin{eqnarray*}
&&(1.0674,-0.1656,1.6253,1.7877,0.3535,2.6909,2.6892,\\
&&\hspace*{5.4pt}1.4624,1.8273,1.6746,1.3133,2.2258,0.9117,3.6832,\\
&&\hspace*{5.4pt}1.3636,1.6139,2.5668,1.5593,1.4044,0.6677,1.7944,0.1638,\\
&&\hspace*{5.4pt}2.2143,3.1236,0.8082,2.7540,-0.0937,0.0590,2.0711,2.3579).
\end{eqnarray*}
Note that $\# (j\dvtx0< |\eta_j |\leq1, j=1,2,\ldots
)=8$, $\# (j\dvtx1< |\eta_j |\leq
2, j=1,2,\ldots)=12$, $\# (j\dvtx2< |\eta_j |\leq
3, j=1,2,\ldots)=8$, and $\# (j\dvtx
|\eta_j |>3, j=1,2,\ldots)=2$.
\end{EX}
\begin{EX}\label{sing.seq.ex6}
We generate the values of $\eta_{j}$, $j=1,\ldots,30$, from the
standard exponential distribution. The
remaining values of $\eta_j$ are 0. The values different from 0 are as follows:
\begin{eqnarray*}
&&(0.0512,1.4647,0.4995,0.7216,0.1151,0.2716,0.7842,\\
&&\hspace*{5.4pt}3.7876,0.1967,0.8103,0.4854,0.2332,0.5814,0.3035,\\
&&\hspace*{5.4pt}1.7357,0.9021,0.0667,0.0867,0.8909,0.1124,2.8491,1.0416,\\
&&\hspace*{8.2pt}0.2068,2.6191,1.9740,1.5957,1.6158,0.5045,1.3012,1.6153).
\end{eqnarray*}
Note that $\# (j\dvtx0<\eta_j\leq1, j=1,2,\ldots)=19$, $\#
(j\dvtx1<\eta_j\leq2, j=1,2,\ldots)=8$,
$\# (j\dvtx2<\eta_j\leq3, j=1,2,\ldots)=2$, and $\# (j\dvtx
\eta_j>3, j=1,2,\ldots)=1$.
\end{EX}
\begin{EX}\label{sing.seq.ex7}
In this example, the values of $\eta_{j}$, $j=1,\ldots,30$, are 2.0
and the rest are zero.
\end{EX}

As expected, the power in each column decreases by increasing $r$
because the number of $\bolds\theta$ with values
different from zero (denoted by $k_{500}^{\mathrm{opt}}$) decreases. When the
$\theta_i$ with the large value such as
3.6832, 3.1236 (in Example \ref{sing.seq.ex4}) and 3.7876 (in Example
\ref{sing.seq.ex6}) is excluded at the
alternative, the large decrement in the power occurs. For each
alternative, the statistic $T_L(15)$ or $T_L(40)$
achieves better power than the order threshold statistics with the
other specified values of the threshold
parameter. This is a consequence of the fact that the number of mean
values that are different from zero never
exceeds 30. Thus, less noise is incorporated in $T_L(k_{500}^{\mathrm{opt}})$
than the other order threshold statistics.
Note that with the chosen value of $\delta_{500}=5.1216$, the hard
threshold statistic uses, on average, 12
observations. Thus, it is rather surprising that the empirical power of
the hard threshold statistic is always
smaller than that of $T_L(15)$. In all three tables, the empirical
power using the approximation $T_L(\widehat
k_{500}^{\mathrm{opt}})\stackrel{\cdot}{\sim} \Phi$ is similar to that of
$T_L(k_{500}^{\mathrm{opt}})$, and always greater than the
empirical powers of the hard threshold and Simes statistics. The
empirical power using the approximation
$T_L(\widehat k_{500}^{\mathrm{opt}})\stackrel{\cdot}{\sim} b\chi^2_\nu$ is
a little bit smaller than that using the normal
approximation, however, it is still greater than the empirical powers
of the hard threshold and Simes
statistics.
%
%
In Table \ref{table7}, for large number of the false null hypotheses the Simes
statistic $T_S$ performs much worse than
the hard threshold statistic, the order threshold statistic, and even
the chi-square statistic $T_L(500)$. In
all three tables, the power of $T_H(5.1216)$ is similar (though
somewhat smaller) to that of $T_L(100)$.
Finally, all order threshold statistics achieved higher power than the
chi-square statistic $T_L(500)$.

\subsection{Choosing $k_n$}\label{choosingkn}

The simulation results and the discussion in the closing paragraph of
Section \ref{ss.sim.sing.seq} suggest that the
power of $T_L(k_n)$ is largest\vspace*{1pt} when $k_n$ equals the number of mean
values different from zero (denoted by
$k_n^{\mathrm{opt}}$). As a data-driven choice of $k_n$, we propose to use the
estimate of $k_n^{\mathrm{opt}}$ suggested by
Storey (\citeyear{S02}, \citeyear{S03}) and Efron et al. (\citeyear{Efronetal01}),
which is
\[
\widehat
k_n^{\mathrm{opt}}(\lambda)=\max\biggl\{\frac{n\mathds{G}_n(\lambda)
-n\lambda-1}{1-\lambda}, \log^{3/2} n \biggr\},
\]
where $\mathds{G}_n$ is the empirical cdf of
$\mathbf{P}^n=(P_1,\ldots,P_n)$, the $P_i$'s are the $p$-values of
the individual hypotheses, and $\lambda$ is the median of the $P_i$'s.
The recommended lower bound $\log^{3/2}n$
of $\widehat k_n^{\mathrm{opt}}(\lambda)$ was found to be preferable in the
simulations we performed. Interestingly,
$\log^{3/2}n$ equals the expected number of observations in hard
thresholding with the recommended threshold
parameter of $\delta_n=2\log(n\log^{-2} n)$.

\section{One-way HANOVA}\label{sec4}

Let the $X_{ij}$, $i=1,\ldots,a$, $j=1,\ldots,n$, be independent
$N(\theta_i,\sigma^2)$, where the $\theta_i$ and $\sigma^2$ are all
unknown. Let $\alpha_i=\theta_i-\overline\theta$ denote the
``effect'' of the $i$th group, and consider testing
$H_0\dvtx\alpha_1=\cdots=\alpha_a=0$ vs. $H_a\dvtx H_0$ is false.
\citet{AP04} show that the asymptotic power of the optimal
invariant ANOVA $F$ test equals its level of significance even when
$\Vert\bolds\alpha\Vert\to\infty$, as $a\to\infty$, with
$\Vert\bolds\alpha\Vert^2=o(\sqrt{a})$. Because the power of the chi-square
statistic (\ref{rel.naive.stat}) has a similar property [\citet{F96}],
an extension of the order thresholding to the one-way HANOVA setting
is expected to result in similar gains in power over the ANOVA $F$
test.

In Section \ref{home.order.test}, we extend the applicability of
order thresholding to the one-way HANOVA context, while Section
\ref{home.sim.result} illustrates the improved power of order
thresholding via simulation. Finally, using the idea of
Storey (\citeyear{S02}, \citeyear{S03}) and the simulation results, we present a recommendation for a
data-driven choice of the order threshold parameter in Section
\ref{HANOVA.choose.ka}.

\subsection{Order thresholding in one-way HANOVA}
\label{home.order.test}

The classical $F$ statistic is given by
%
%
\begin{equation}\label{paper1-6.3}
F_a=\frac{\mathrm{MST}}{\mathrm{MSE}},
\end{equation}
where
\[
\mathrm{MST}=\frac{1}{a-1}\sum_{i=1}^a n (\overline{X}_{i\cdot}-\overline
{X}_{\cdot\cdot} )^2,\qquad
\mathrm{MSE}=\frac{1}{N-a}\sum_{i=1}^a\sum_{j=1}^{n} (X_{ij}-\overline
{X}_{i\cdot} )^2
\]
with $\overline{X}_{i\cdot}=n^{-1}\sum_{j=1}^{n}X_{ij}$,
$\overline{X}_{\cdot\cdot}=N^{-1}\sum_{i=1}^a\sum_{j=1}^{n}X_{ij}$, and
$N=an$. Note that
%
%
\begin{equation}\label{paper1-6.4}
(a-1)F_a=\sum_{i=1}^a \biggl(\frac{\sqrt{n} (\overline{X}_{i\cdot}
-\overline{X}_{\cdot\cdot} )}
{\sqrt{\mathrm{MSE}}} \biggr)^2
\end{equation}
differs from the chi-square statistic (\ref{rel.naive.stat}) only in
that the random variables which are being
summed are not independent, and their distribution is not $\chi^2_1$. Set
\[
\widetilde{Z}_i=\frac{\sqrt{n} (\overline{X}_{i\cdot} -\overline
{X}_{\cdot\cdot} )} {\sqrt{\mathrm{MSE}}}, \qquad
\widehat Z_i=\frac{\sqrt{n} (\overline{X}_{i\cdot} -\overline
{X}_{\cdot\cdot} )} {\sigma}.
\]
Thus,
\[
\widetilde{Z}_i=s\widehat Z_i\qquad \mbox{where } s=\frac{\sigma
}{\sqrt{\mathrm{MSE}}}.
\]
Threshold versions of (\ref{paper1-6.4}) are of the form
%
%
\begin{equation}\label{rel.F.thres}
\widehat T_L(k_a)=\sum_{i=1}^ac_{ia}\widetilde{Z}_{i,a}^2 =s^2\sum
_{i=1}^ac_{ia}\widehat{Z}_{i,a}^2,
\end{equation}
where $\widehat{Z}_{1,a}^2<\cdots<\widehat{Z}_{a,a}^2$ are the
ordered $\widehat{Z}_{i}^2$'s,
$\widetilde{Z}_{i,a}=s\widehat{Z}_{i,a}$, $c_{ia}=I(i>a-k_a)$, and
$k_a$ is the order threshold parameter. For
suitable centering and scaling constants, $\widehat\mu_a(k_a)$ and
$\widehat\sigma_a(k_a)$, the asymptotic
theory of $\widehat T_L(k_a)$ will use the decomposition
%
%
\begin{eqnarray}\label{rel.decomp}
\frac{\widehat T_L(k_a)-a\widehat\mu_a(k_a)}
{\sqrt{a}\widehat\sigma_a(k_a)}&=&s^2\frac{1}{\sqrt{a}\widehat
\sigma_a(k_a)}
\Biggl(\sum_{i=1}^ac_{ia}\widehat{Z}_{i,a}^2-a\widehat\mu
_a(k_a) \Biggr)\nonumber\\[-8pt]\\[-8pt]
&&{}+\frac{\sqrt{a}}{\widehat\sigma_a(k_a)}\widehat\mu_a(k_a)
(s^2-1 ).\nonumber
\end{eqnarray}
The two components in (\ref{rel.decomp}) are independent, so it
suffices to show the asymptotic normality of
each one separately. To deal with the first component, let $\theta_0$
denote the common value of the $\theta_i$
under $H_0$ and write
%
%
\begin{equation}\label{rel.st}
\widehat{Z}_i=Z_i+\frac{t}{\sqrt{a}}\qquad \mbox{where }
Z_i=\frac{\sqrt{n} (\overline{X}_{i\cdot}-\theta_0
)}{\sigma} \mbox{ and }
t=-\frac{\sqrt{N} (\overline{X}_{\cdot\cdot}-\theta_0
)}{\sigma}.\hspace*{-27pt}
\end{equation}
Our approach for obtaining the asymptotic distribution of $\sum
_{i=1}^ac_{ia}\widehat Z_{i,a}^2$ is to first
derive its asymptotic distribution treating the $t$ in (\ref{rel.st})
as fixed, and then to show that the
convergence is uniform over all values of $t$ bounded by any positive
constant $M$. By Slutsky's theorem,\vspace*{1pt} the
asymptotic distribution of the first component of (\ref{rel.decomp})
is the same as that of
$\sum_{i=1}^ac_{ia}\widehat Z_{i,a}^2$. The asymptotic distribution of
the second component of
(\ref{rel.decomp}) is easily derived since
\[
\sqrt{a}(s^2-1)\stackrel{d}{\to}N \biggl(0,\frac{2}{n-1}
\biggr)\qquad
\mbox{as } a\to\infty,
\]
and, as it will be shown in lemmas described in Section \ref{anova.order},
$\widehat\mu_a(k_a)/\widehat\sigma_a(k_a)\to\mu_r/\sigma_r$,
provided that $k_a/a\to r$ for some $0\leq r\leq
1$ and $k_a\to\infty$, as $a\to\infty$.

\subsubsection{Asymptotic distribution when $t$ is fixed}

When $t$ is fixed, we set
%
%
\begin{equation}\label{paper1-6.10}
Z_{t,i}=Z_i+\frac{t}{\sqrt{a}},\qquad i=1,\ldots,a,\qquad T_L^{t}(k_a)=
\sum_{i=1}^ac_{ia}Z_{t,(i)}^2,
\end{equation}
where $Z_{t,(1)}^2<\cdots<Z_{t,(a)}^2$ are the order statistics of
$Z_{t,1}^2,\ldots,Z_{t,a}^2$. [Note that for
$t$ as defined in (\ref{rel.st}), $Z_{t,i}$ becomes $\widehat Z_i$.]
It follows that the $Z_{t,i}^2$ are
independent $\chi^2_1 (t^2/a )$ so that their density and
cumulative distribution functions are given
by
\[
g_{a,t}(y)=\frac{e^{-{1}/{2} (y+{t^2}/{a}
)}y^{-1/2}}{2^{1/2}}
\sum_{k=0}^{\infty}\frac{({t^2}/{a}y)^k}{2^{2k}k!\Gamma
(k+1/2 )},\qquad y>0
\]
and
\[
G_{a,t}(y)=\int_0^y
g_{a,t}(u)\,du=\sum_{k=0}^{\infty}e^{-{t^2}/({2a})}\frac
{1}{2^k\cdot k!} \biggl(\frac{t^2}{a} \biggr)^k
G_{2k+1}(y),\qquad y>0,
\]
respectively, where
\[
G_k(y)=\frac{1}{2^{k/2}\Gamma({k/2})}\int
_0^yu^{k/2-1}e^{-u/2}\,du, \qquad y>0
\]
is the cumulative distribution function of $\chi_k^2(0)$. Let
%
%
\begin{equation}\label{eq:3}
\mu_a^t(k_a)=\frac{1}{a}\sum_{i=1}^a
c_{ia}G_{a,t}^{-1}(1-e^{-\tilde{\nu}_{ia}}) \quad\mbox{and}\quad
(\sigma_a^t(k_a))^2=\frac{1}{a}\sum_{i=1}^a (\alpha
_{ia}^t(k_a) )^2,\hspace*{-25pt}
\end{equation}
where $\alpha_{ia}^t(k_a)=\frac{1}{a-i+1}\sum_{j=i}^ac_{ja}\frac{
e^{-\tilde{\nu}_{ja}}}
{g_{a,t} (G_{a,t}^{-1}(1-e^{-\tilde{\nu}_{ja}}) )}$ and
$\tilde{\nu}_{ia}=\sum_{j=1}^i\frac{1}{a-j+1}$. With this notation
we have the following lemma.
\begin{LEM}[{[\citet{CGJ67}]}]\label{anova.lem.4.1}
Let $T_L^t(k_a)$ and $\mu_a^t(k_a)$ be as defined in (\ref
{paper1-6.10}) and (\ref{eq:3}), respectively. Let
$V_1,\ldots,V_a$ be i.i.d. from $\operatorname{Exp}(1)$ random variables and let
$V_{1,a}<\cdots<V_{a,a}$ be the corresponding
order statistics. Then $a^{-1}T_L^t(k_a)$ can be decomposed as
\[
a^{-1}T_L^{t}(k_a)\stackrel{d}{=}\mu_a^t(k_a)+Q_a^t(k_a)+R_a^t(k_a),
\]
where
%
%
\begin{equation}\label{rel.def.Q}
Q_a^t(k_a)=\frac{1}{a}\sum_{i=1}^a\alpha_{ia}^t(k_a)(V_i-1)
\end{equation}
and
\[
R_a^t(k_a)=\frac{1}{a}\sum_{i=1}^a
c_{ia} \biggl\{ \bigl(G_{a,t}^{-1} (1-e^{-V_{i,a}} )-G_{a,t}^{-1}
(1-e^{-{\tilde{\nu}}_{ia}} ) \bigr)-\frac
{(V_{i,a}-\tilde{\nu}_{ia})
e^{-\tilde{\nu}_{ia}}}{g_{a,t} (G_{a,t}^{-1}
(1-e^{-\tilde{\nu}_{ia}} ) )} \biggr\}
\]
with
$\alpha_{ia}^t(k_a)=\frac{1}{a-i+1}\sum_{j=i}^ac_{ja}\frac{
e^{-\tilde{\nu}_{ja}}}
{g_{a,t} (G_{a,t}^{-1}(1-e^{-\tilde{\nu}_{ja}}) )}$ and
$\tilde{\nu}_{ia}=\sum_{j=1}^i\frac{1}{a-j+1}$.
\end{LEM}
\begin{Theorem}\label{anova.thm.4.1}
For any fixed value of $t$, let $Z_{t,i}^2$, $i=1,\ldots,a$, be a
sequence of i.i.d. random variables having the noncentral chi-squared
distribution with 1 degree of freedom and noncentrality parameter
$t^2/a$. Let $k_a, a\ge1$, be any sequence of integers which
satisfies $k_a\to\infty$, as $a\to\infty$, and $k_a\le a$. Let
$\mu_a^t(k_a)$ and $ (\sigma_a^t(k_a) )^2$ be as in
(\ref{eq:3}) with $c_{ia}=I(i>a-k_a)$, and let $T_L^t(k_a)$ be given
in (\ref{paper1-6.10}). Then we have 
%
%
\begin{equation}\label{paper1-6.15}
{T_L^{t}}^*(k_a)=\frac{T_L^{t}(k_a)-a\mu_a^t(k_a)} {\sqrt{a}\sigma
_a^t(k_a)}\stackrel{d}{\rightarrow}N(0,1)\qquad
\mbox{as } a\to\infty.
\end{equation}
%
\end{Theorem}

\subsubsection{Uniformity of the convergence in distribution}

This subsection shows that the distribution function of (\ref
{paper1-6.15}) converges to the standard normal
distribution
uniformly on $|t|<M$. 
%
\begin{LEM}\label{anova.lem.4.2}
Consider the setting of Theorem \ref{anova.thm.4.1}. Let $H_{a,t}$ be
the distribution function of
$\sqrt{a}Q_a^t(k_a)/\sigma_a^t(k_a)$, where $Q_a^t(k_a)$ is given in
(\ref{rel.def.Q}), and let $\Phi$ be the
standard normal distribution function. Then, for any $M>0$,
\[
{\mathop{\sup_{-M< t< M}}_{-\infty< x<\infty}}
|H_{a,t}(x)-\Phi(x) |\rightarrow0\qquad\mbox{as }
a\rightarrow\infty.
\]
\end{LEM}
\begin{LEM}\label{anova.lem.4.3}
Consider the setting of Theorem \ref{anova.thm.4.1}, and let
$R_a^t(k_a)$ be as given in Lemma
\ref{anova.lem.4.1}. Then, for any $M>0$,
\[
\sup_{-M< t< M}
\biggl|\frac{\sqrt{a}R_a^t(k_a)}{\sigma_a^t(k_a)} \biggr|\stackrel
{p}{\rightarrow}0\qquad
\mbox{as } a\rightarrow\infty.
\]
\end{LEM}
\begin{LEM}\label{anova.thm.4.2}
Consider the setting of Theorem \ref{anova.thm.4.1}. Let $F_{a,t}$ be
the distribution function of
${T_L^{t}}^*(k_a)$ given in (\ref{paper1-6.15}) and let $\Phi$ be the
standard normal distribution function.
Then, for any $M>0$,
\[
{\mathop{\sup_{-M< t< M}}_{-\infty< x<\infty}}
|F_{a,t}(x)-\Phi(x) | \rightarrow0\qquad \mbox{as }
a\rightarrow\infty.
\]
\end{LEM}
\begin{Theorem}\label{anova.t.random} Let $k_a, a\ge1$, be any
sequence of integers which satisfies $k_a\to\infty$, as
$a\to\infty$, and $k_a\le a$. For $t$ as defined in (\ref{rel.st}),
let $\widehat{Z}_i$, $\widehat\mu_a(k_a)$ and
$ (\widehat\sigma_a(k_a) )^2$ be as in (\ref{paper1-6.10}),
(\ref{eq:3}), respectively. Then we have
%
%
\begin{equation}
\widehat T_L^*(k_a)=\frac{\sum_{i=1}^ac_{ia}\widehat
Z_{i,a}^2-a\widehat{\mu}_a(k_a)}{\sqrt{a}\widehat{\sigma
}_a(k_a)}\stackrel{d}{\to} N(0,1)\qquad \mbox{as }
a\to\infty,
\end{equation}
where $\widehat{Z}_{1,a}^2<\cdots<\widehat{Z}_{a,a}^2$ are the
ordered $\widehat{Z}_{i}^2$'s and $c_{ia}=I(i>a-k_a)$.
\end{Theorem}

\subsubsection{Asymptotic normality of the order threshold
statistics}\label{anova.order}

In this subsection, it is first shown that $\mu_a^t(k_a)$ and
$\sigma_a^t(k_a)$ converge to $\mu_a^0(k_a)$ and $\sigma_a^0(k_a)$,
respectively, uniformly on $|t|<M$. This fact is then used in
Theorem \ref{anova.order.thres} for obtaining the asymptotic
normality of the order threshold statistic given in
(\ref{rel.F.thres}).
\begin{LEM}\label{ratio.std}
Let $k_a, a\ge1$, be any sequence of integers which satisfies
$k_a\to\infty$, as $a\to\infty$, and $k_a\le a$. Let
$ (\sigma_a^t(k_a) )^2$ and $(\sigma_a^0(k_a))^2$ be as in
(\ref{eq:3}) with any $t$ and fixed value of $t=0$, respectively.
Then, for any $M>0$,
\[
\sup_{-M< t< M}
\biggl|\frac{\sigma_a^t(k_a)}{\sigma_a^0(k_a)}-1
\biggr|\rightarrow0\qquad
\mbox{as } a\rightarrow\infty.
\]
\end{LEM}
\begin{LEM}\label{diff.mu}
Let $k_a, a\ge1$, be any sequence of integers which satisfies
$k_a\to\infty$, as $a\to\infty$, and $k_a\le a$. Let $\mu_a^t(k_a)$,
$\mu_a^0(k_a)$ and $(\sigma_a^0(k_a))^2$ be as in (\ref{eq:3}) with
any $t$, fixed value of $t=0$, respectively. Then, for any $M>0$,
\[
\sup_{-M< t< M}
\biggl|\frac{\sqrt{a} (\mu_a^t(k_a)-\mu_a^0(k_a)
)}{\sigma_a^0(k_a)} \biggr| \rightarrow0
\qquad\mbox{as } a\rightarrow\infty.
\]
\end{LEM}
\begin{LEM}\label{second.lem}
Let $\mu_a^0(k_a)$ and $(\sigma_a^0(k_a))^2$ be as in (\ref{eq:3})
with the fixed value of $t=0$. Then, provided that $k_a/a\to r$ for
some $0\leq r\leq1$ and $k_a\to\infty$, as $a\to\infty$, we have
\[
\mu_a^0(k_a)\to\mu_r \quad\mbox{and}\quad (\sigma_a^0(k_a))^2\to
\sigma_r^2\qquad \mbox{as } a\to\infty,
\]
where
\[
\mu_r=\int_0^1I(t>1-r)G_{a,0}^{-1}(t)\,dt
\]
and
\[
\sigma_r^2=\int_0^1\int_0^1I(t>1-r)I(s>1-r)\bigl(\min(t,s)-ts\bigr)
\,dG_{a,0}^{-1}(t)\,dG_{a,0}^{-1}(s).
\]
\end{LEM}
\begin{RE*}
If $r=1$, then
\[
\frac{\mu_a^0(k_a)}{\sigma_a^0(k_a)} \to\frac{1}{\sqrt{2}}\qquad
\mbox{as } a\to\infty.
\]
\end{RE*}

From Theorem \ref{anova.t.random} and lemmas described earlier in this
subsection, we can obtain the following theorem.
%
%
\begin{Theorem}\label{anova.order.thres}
Let $\mu_a^0(k_a)$, $(\sigma_a^0(k_a))^2$, $\mu_r$, and $\sigma
_r^2$ be as in Lemma \ref{second.lem}, and let
$\widehat T_L(k_a)$ be given in (\ref{rel.F.thres}). Then, provided
that $k_a/a\to r$ for some $0\leq r\leq1$
and $k_a\to\infty$, as $a\to\infty$, we have
%
%
\begin{equation}\label{paper1-6.21}
\widetilde T_L(k_a)=\frac{\widehat T_L(k_a) -a\mu_a^0(k_a)}
{\sqrt{a}\sigma_a^0(k_a)}\stackrel{d}{\rightarrow} N
\biggl(0,1+\frac{2\mu_r^2}{\sigma_r^2(n-1)} \biggr)\qquad
\mbox{as } a\to\infty.\hspace*{-27pt}
\end{equation}
\end{Theorem}

\subsection{Simulations}\label{home.sim.result}

In this subsection, we compare the performance of the classical $F$
statistic, given in (\ref{paper1-6.3}), and the order threshold
statistics $\widetilde T_L(k_a)$, given in (\ref{paper1-6.21}).

We remark that \citet{FL98} applied the thresholding
methodology to the problem of comparing $I$ curves with data arising
from the model $X_{ij}(t)=f_i(t)+ \varepsilon_{ij}(t), t=1,\ldots,T,
j=1,\ldots,n_i, i=1,\ldots,I$. Their asymptotic theory pertains to
the case where the number of curves which are compared, $I$, remains
fixed, while $T$ and the sample sizes $n_i$ tend to infinity. This
problem is fundamentally different from that considered here, and
their procedure is not a competitor to ours.

\begin{figure}

\includegraphics{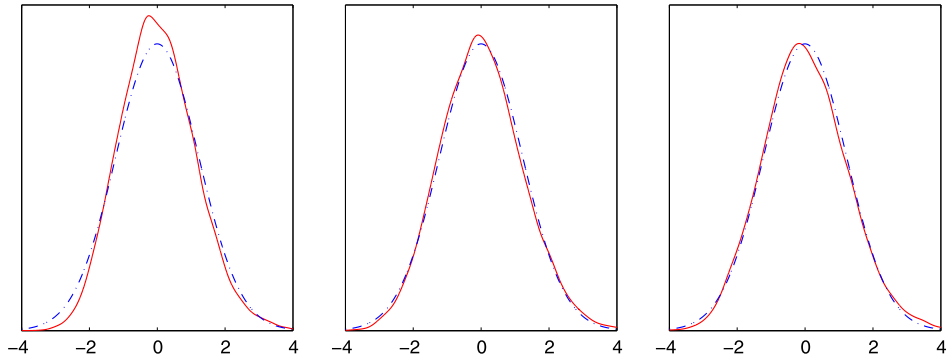}

\caption{Estimated densities of $\widetilde T_L(k_{500})$ for $a=500$,
$n=3$, and $k_{500}= 22$, 105 and 229.}
\label{fig.order thres anova_pp}
\end{figure}

Figure \ref{fig.order thres anova_pp} presents the estimated
densities of $\widetilde T_L(k_{500})$ (solid line) and the density
of the limiting normal distribution (dash-dot line). The estimated
densities are based on 20,000 simulated values, using $a=500$ and
$n=3$, when the threshold parameter $k_{500}$ takes the values of
$[a^{1/2}]=22$, $[a^{3/4}]=105$, and $[a^{7/8}]=229$. It can be seen
that the approximation is quite good especially for $k_{500}=105$
and 229. Similar figures (not shown here) with different values of
$n$ suggest that the rate of convergence of the order threshold
statistic to its limiting distribution is mainly driven by $a$, not
$n$.

%
%

%
%
\begin{table}[b]
\tabcolsep=0pt
\caption{Type \textup{I} errors of order threshold statistics, $\widetilde
T_L(k_a)$, for different values of the
threshold parameter}\label{paper1-order-anova2}
%
\begin{tabular*}{\tablewidth}{@{\extracolsep{\fill}}lcccccccc@{}}
\hline
&$\bolds{[\log^{1/2} a]}$& $\bolds{[\log a]}$ & $\bolds{[\log^{3/2} a]}$
& $\bolds{[a^{1/2}]}$ & $\bolds{[a^{2/3}]}$ & $\bolds{[a^{3/4}]}$ &
$\bolds{[a^{7/8}]}$ & $\bolds{a}$\\
\hline
$a=50$ and $n=3$ & 0.0522& 0.0551& 0.0601& 0.0601& 0.0623& 0.0635&
0.0637& 0.0669\\
$a=50$ and $n=5$ & 0.0551& 0.0583& 0.0591& 0.0591& 0.0588& 0.0600&
0.0612& 0.0619\\
$a=100$ and $n=3$ & 0.0506& 0.0521& 0.0561& 0.0563& 0.0594& 0.0607&
0.0617& 0.0634\\
$a=100$ and $n=5$ & 0.0539& 0.0541& 0.0541& 0.0549& 0.0571& 0.0578&
0.0596& 0.0604\\
$a=200$ and $n=3$ & 0.0436& 0.0440& 0.0490& 0.0497& 0.0552& 0.0571&
0.0601& 0.0597\\
$a=200$ and $n=5$ & 0.0548& 0.0520& 0.0505& 0.0504& 0.0515& 0.0529&
0.0542& 0.0549\\
$a=500$ and $n=3$ & 0.0436& 0.0437& 0.0452& 0.0466& 0.0515& 0.0558&
0.0593& 0.0589\\
$a=500$ and $n=5$ & 0.0533& 0.0492& 0.0474& 0.0481& 0.0510& 0.0518&
0.0532& 0.0534\\
$a=1000$ and $n=3$ & 0.0427& 0.0403& 0.0405& 0.0411& 0.0475& 0.0513&
0.0548& 0.0557\\
$a=1000$ and $n=5$ & 0.0517& 0.0486& 0.0459& 0.0453& 0.0466& 0.0484&
0.0507& 0.0521\\
\hline
\end{tabular*}
\end{table}

The results reported in Table \ref{paper1-order-anova2} are based on 20,000 simulation runs. As
expected, the distributions of
$\widetilde T_L(k_a)$ converge to the normal distribution function and
the achieved alpha levels are close to
the true value of 0.05. Thus, the asymptotic theory of the order
threshold statistics provides a good
approximation. More exactly, when the number of groups are larger than
200, all order threshold statistics are
robust for the 0.05 significance level. In particular, the achieved
alpha level of $\widetilde T_L(k_{1000})$
is 0.0507 when $a=1000$, $n=5$, and $k_{1000}=[a^{7/8}]=421$.

From now, we compare the empirical power of $\widetilde T_L(k_{1000})$
using several values of the threshold
parameter with that of the classical $F$ statistic. The simulations use
samples of size $a=1000$ and $n=5$
generated from the normal distribution with variance 1. The threshold
parameter $k_{1000}$ is 20, 50, 100, 250,
500, and 1000. All results are based on 20,000 simulation runs. The
alternatives here have 20 of the 1000
$\theta_i$ values different from zero. In particular, we consider the
following sequence of alternatives
indexed by $r$:
\[
H_r\dvtx\theta_j=\eta_{j+r-1} \qquad\mbox{for } j=1,\ldots,1000,
r=1,\ldots,20,
\]
where $\eta_j$, $j=1,2,\ldots,$ is a given sequence. The following
are examples with different values of
$\bolds\eta$.
\begin{EX}\label{anova.seq.ex1}
We generate the values of $\eta_{j}$, $j=1,\ldots,20$, from\break
$\operatorname{Uniform}(-2,2)$. The remaining values of $\eta_j$
are 0. The values different from 0 are as follows:
\begin{eqnarray*}
&&(1.8005, -1.0754, 0.4274, -0.0561, 1.5652, 1.0484, \\
&&\hspace*{5.4pt}-0.1741,-1.9260, 1.2856, -0.2212, 0.4617, 1.1677,1.6873, \\
&&\hspace*{8.4pt} 0.9528, -1.2949, -0.3772, 1.7419, 1.6676,
-0.3589, 1.5746).
\end{eqnarray*}
Note that $\# (j\dvtx0< |\eta_j |\leq1, j=1,2,\ldots
)=8$ and $\# (j\dvtx|\eta_j |>1,
j=1,2,\ldots)=12$.
\end{EX}
\begin{EX}\label{anova.seq.ex2}
We generate the values of $\eta_{j}$, $j=1,\ldots,20$, from
$\operatorname{Exp}(0.7)$. The remaining values of $\eta_j$ are 0.
The values different from 0 are as follows:
\begin{eqnarray*}
&&(1.0949, 0.5511, 1.7587, 0.1128, 0.4033, 0.7991, 0.6868,\\
&&\hspace*{5.4pt} 0.0993,0.6919, 1.8255, 1.1272,2.1041, 0.3975,\\
&&\hspace*{6.8pt} 1.4730, 0.4549, 1.5015, 0.1830, 0.6865, 0.1360, 2.1458).
\end{eqnarray*}
Note that $\# (j\dvtx0<\eta_j\leq1, j=1,2,\ldots)=12$, $\#
(j\dvtx1<\eta_j\leq2, j=1,2,\ldots)=6$
and $\# (j\dvtx\eta_j> 2, j=1,2,\ldots)=2$.
\end{EX}

%
%
\begin{table}
\caption{Power calculations in
Example \protect\ref{anova.seq.ex1}}\label{paper1-anova-seq-t1}
\begin{tabular*}{\tablewidth}{@{\extracolsep{\fill}}lcccccccc@{}}
\hline
&$\bolds{k_{1000}^{\mathrm{opt}}}$&$\bolds{F_{1000}}$ &$\bolds{\widetilde T_L(20)}$
&$\bolds{\widetilde T_L(50)}$ & $\bolds{\widetilde
T_L(100)}$ & $\bolds{\widetilde T_L(250)}$ & $\bolds{\widetilde T_L(500)}$ &
$\bolds{\widetilde T_L(1000)}$\\
\hline
$H_{1}$& 20& 0.8612 & 0.9992 & 0.9975 & 0.9877 & 0.9482 & 0.8923 &
0.8682\\
$H_{2}$& 19& 0.7887 & 0.9963 & 0.9889 & 0.9685 & 0.9000 & 0.8270 &
0.7978\\
$H_{3}$& 18& 0.7561 & 0.9957 & 0.9878 & 0.9623 & 0.8762 & 0.7971 &
0.7658\\
$H_{4}$& 17& 0.7505 & 0.9952 & 0.9848 & 0.9588 & 0.8743 & 0.7924 &
0.7601\\
$H_{5}$& 16& 0.7541 & 0.9949 & 0.9841 & 0.9591 & 0.8801 & 0.7944 &
0.7633\\
$H_{6}$& 15& 0.6785 & 0.9901 & 0.9712 & 0.9275 & 0.8175 & 0.7238 &
0.6891\\
$H_{7}$& 14& 0.6434 & 0.9859 & 0.9634 & 0.9116 & 0.7856 & 0.6887 &
0.6563\\
$H_{8}$& 13& 0.6432 & 0.9855 & 0.9623 & 0.9100 & 0.7876 & 0.6905 &
0.6547\\
$H_{9}$& 12& 0.5091 & 0.9422 & 0.8861 & 0.8008 & 0.6505 & 0.5518 &
0.5193\\
$H_{10}$& 11& 0.4434 & 0.9191 & 0.8399 & 0.7351 & 0.5794 & 0.4868 &
0.4553\\
$H_{11}$& 10& 0.4444 & 0.9191 & 0.8399 & 0.7355 & 0.5742 & 0.4855 &
0.4561\\
$H_{12}$& \phantom{0}9& 0.4448 & 0.9230 & 0.8414 & 0.7333 & 0.5760 & 0.4847 &
0.4562\\
$H_{13}$& \phantom{0}8& 0.3896 & 0.8894 & 0.7869 & 0.6756 & 0.5132 & 0.4264 &
0.4007\\
$H_{14}$& \phantom{0}7& 0.2887 & 0.7710 & 0.6364 & 0.5169 & 0.3835 & 0.3185 &
0.2989\\
$H_{15}$& \phantom{0}6& 0.2615 & 0.7437 & 0.6051 & 0.4866 & 0.3537 & 0.2903 &
0.2724\\
$H_{16}$& \phantom{0}5& 0.2095 & 0.6603 & 0.5037 & 0.3878 & 0.2803 & 0.2321 &
0.2187\\
$H_{17}$& \phantom{0}4& 0.2089 & 0.6560 & 0.5002 & 0.3869 & 0.2742 & 0.2319 &
0.2169\\
$H_{18}$& \phantom{0}3& 0.1356 & 0.4002 & 0.2874 & 0.2250 & 0.1686 & 0.1482 &
0.1421\\
$H_{19}$& \phantom{0}2& 0.0816 & 0.1736 & 0.1287 & 0.1106 & 0.0943 & 0.0884 &
0.0867\\
$H_{20}$& \phantom{0}1& 0.0812 & 0.1743 & 0.1277 & 0.1095 & 0.0934 & 0.0880 &
0.0862\\
\hline
\end{tabular*}
%
\end{table}

%
%
\begin{table}
\caption{Power calculations in
Example \protect\ref{anova.seq.ex2}}\label{paper1-anova-seq-t2}
\begin{tabular*}{\tablewidth}{@{\extracolsep{\fill}}lcccccccc@{}}
\hline
& $\bolds{k_{1000}^{\mathrm{opt}}}$ & $\bolds{F_{1000}}$ & $\bolds{\widetilde T_L(20)}$
& $\bolds{\widetilde T_L(50)}$ & $\bolds{\widetilde T_L(100)}$
& $\bolds{\widetilde T_L(250)}$ & $\bolds{\widetilde T_L(500)}$
& $\bolds{\widetilde T_L(1000)}$\\
\hline
$H_{1}$& 20& 0.7680& 0.9978& 0.9886& 0.9657& 0.8877& 0.8089& 0.7769\\
$H_{2}$& 19& 0.7275& 0.9968& 0.9861& 0.9550& 0.8603& 0.7732& 0.7366\\
$H_{3}$& 18& 0.7241& 0.9960& 0.9842& 0.9533& 0.8563& 0.7669& 0.7330\\
$H_{4}$& 17& 0.6278& 0.9893& 0.9640& 0.9048& 0.7740& 0.6731& 0.6394\\
$H_{5}$& 16& 0.6253& 0.9886& 0.9624& 0.9052& 0.7702& 0.6730& 0.6373\\
$H_{6}$& 15& 0.6188& 0.9892& 0.9624& 0.9031& 0.7681& 0.6667& 0.6306\\
$H_{7}$& 14& 0.6011& 0.9872& 0.9577& 0.8891& 0.7464& 0.6462& 0.6119\\
$H_{8}$& 13& 0.5871& 0.9872& 0.9519& 0.8829& 0.7369& 0.6337& 0.5982\\
$H_{9}$& 12& 0.5831& 0.9870& 0.9530& 0.8819& 0.7406& 0.6342& 0.5962\\
$H_{10}$& 11& 0.5614& 0.9849& 0.9467& 0.8730& 0.7151& 0.6097& 0.5750\\
$H_{11}$& 10& 0.4476& 0.9526& 0.8704& 0.7600& 0.5872& 0.4900& 0.4598\\
$H_{12}$& \phantom{0}9& 0.4009& 0.9411& 0.8435& 0.7224& 0.5399& 0.4405& 0.4121\\
$H_{13}$& \phantom{0}8& 0.2521& 0.7461& 0.5879& 0.4612& 0.3297& 0.2770& 0.2612\\
$H_{14}$& \phantom{0}7& 0.2495& 0.7446& 0.5843& 0.4573& 0.3319& 0.2742& 0.2597\\
$H_{15}$& \phantom{0}6& 0.1831& 0.6204& 0.4465& 0.3361& 0.2419& 0.2026& 0.1913\\
$H_{16}$& \phantom{0}5& 0.1820& 0.6119& 0.4411& 0.3383& 0.2407& 0.2014& 0.1898\\
$H_{17}$& \phantom{0}4& 0.1283& 0.4346& 0.2941& 0.2197& 0.1613& 0.1412& 0.1356\\
$H_{18}$& \phantom{0}3& 0.1296& 0.4389& 0.2959& 0.2195& 0.1654& 0.1434& 0.1363\\
$H_{19}$& \phantom{0}2& 0.1195& 0.4202& 0.2793& 0.2084& 0.1515& 0.1308& 0.1258\\
$H_{20}$& \phantom{0}1& 0.1176& 0.4207& 0.2763& 0.2041& 0.1532& 0.1296& 0.1238\\
\hline
\end{tabular*}
%
\end{table}

As expected, the power in each column decreases as $r$ increases and
$\widetilde T_L(20)$ has the highest
power. Since the\vspace*{1pt} number of $\theta_i$'s that are different from zero
does not exceed 20, $\widetilde T_L(20)$
minimizes the accumulation of noise, compared to the other order
threshold statistics. For each alternative,
the largest power differences between $F_{1000}$ and $\widetilde
T_L(20)$ are about 0.5 (alternative $H_{13}$
in Table \ref{paper1-anova-seq-t1}) and 0.54 (alternative $H_{12}$ in Table \ref{paper1-anova-seq-t2}). In both
tables, the power of $\widetilde T_L(1000)$ is
similar to that of $F_{1000}$ because $\widetilde T_L(1000)$ is a
standardized version of $F_{1000}$. Finally,
all order threshold statistics achieved higher power than the classical
$F$ statistic $F_{1000}$.

\subsection{Choosing $k_a$}\label{HANOVA.choose.ka}

The simulation results and the discussion in the closing paragraph of
Section \ref{home.sim.result} suggest that
choosing $k_a$ equal to the number of groups with nonzero effects,
$k_a^{\mathrm{opt}}$, maximizes the power. Our
recommendation for the choice of the threshold parameter is based again
on the idea of Storey (\citeyear{S02}, \citeyear{S03}) for
enhancing the power of Simes statistic for testing the constructed set
of hypothesis testing problems
$H_{0}^{(i)}\dvtx\theta_i=\overline{X}_{\cdot\cdot}$, $i=1,\ldots,a$,
where $\overline{X}_{\cdot\cdot}$ is the
overall sample mean. The $p$-value for each hypothesis is approximated by
\[
P_i=2\bigl(1-\Phi(|Z_i|)\bigr),\qquad i=1,\ldots,a,
\]
with
\[
Z_i=\frac{\overline{X}_{i\cdot}-\overline{X}_{\cdot\cdot}} {\sqrt
{S_p^2/n}},\qquad i=1,\ldots,a,
\]
where $\overline{X}_{i\cdot}$ is the sample mean from the $i$th group
and $S_p^2$ is the pooled sample
variance. The power-enhanced version of the Simes statistic
\[
T_S=\min_{1\leq i\leq a} \bigl\{aP_{(i)}/i \bigr\}
\]
rejects the global null hypothesis if $T_S<\alpha/(1-\widehat
k_a^{\mathrm{opt}}/a)$, with
%
%
\begin{equation}\label{est.ka}
\widehat k_a^{\mathrm{opt}}(\lambda)= \max\biggl\{\frac{a\mathds
{G}_a(\lambda)-a\lambda-1}{1-\lambda}, \log^{3/2}
a \biggr\},
\end{equation}
where $\mathds{G}_a$ is the empirical cdf of $\mathbf{P}^a=
(P_1,\ldots,P_a)$, $P_{(1)}<\cdots<P_{(a)}$ are the
ordered $P_i$'s, and $\lambda$ is the median of the $P_i$'s.

The simulation results shown in Table \ref{table11} suggest that the power of
$\widetilde T_L(\widehat k_{1000}^{\mathrm{opt}})$ is
similar to that of $\widetilde T_L(k_{1000}^{\mathrm{opt}})$. These results are
based on 2000 simulation runs; the type
I error rate of $\widetilde T_L(\widehat k_{1000}^{\mathrm{opt}})$ was 0.048.



%
%
\begin{table}
\caption{Power calculations in Example
\protect\ref{anova.seq.ex1}}\label{table11}
\begin{tabular*}{\tablewidth}{@{\extracolsep{\fill}}lcccccccccc@{}}
\hline
&$\bolds{H_{1}}$ & $\bolds{H_{2}}$ & $\bolds{H_{3}}$ & $\bolds{H_{4}}$
& $\bolds{H_{5}}$ & $\bolds{H_{6}}$ & $\bolds{H_{7}}$ &
$\bolds{H_{8}}$ & $\bolds{H_{9}}$ & $\bolds{H_{10}}$\\
\hline
$k_{1000}^{\mathrm{opt}}$& 20 & 19& 18& 17& 16& 15& 14& 13& 12& 11\\
$\widetilde T_L(\widehat k_{1000}^{\mathrm{opt}})$& 1.000& 0.996& 0.993& 0.994&
0.995& 0.987& 0.980& 0.981&
0.938& 0.904\\
\hline
& $\bolds{H_{11}}$ & $\bolds{H_{12}}$ & $\bolds{H_{13}}$ & $\bolds{H_{14}}$
& $\bolds{H_{15}}$ & $\bolds{H_{16}}$ & $\bolds{H_{17}}$ & $\bolds{H_{18}}$
& $\bolds{H_{19}}$ & $\bolds{H_{20}}$\\
\hline
$k_{1000}^{\mathrm{opt}}$& 10 & 9& 8& 7& 6& 5& 4& 3& 2& 1\\
$\widetilde T_L(\widehat k_{1000}^{\mathrm{opt}})$& 0.911& 0.900& 0.883& 0.774&
0.722& 0.674& 0.661& 0.398& 0.182 & 0.175\\
\hline
\end{tabular*}
%
\end{table}


\section{Discussion}\label{sec5}

The asymptotic theory of test statistics based on hard and soft
thresholding pertain the normal distribution
and require the threshold parameter to tend to infinity at a strictly
prescribed rate. This second feature
results in potentially compromised power of the hard threshold statistic.

Order thresholding, a new thresholding method based on order
statistics, is proposed. The asymptotic theory,
developed under the normal distribution in this paper, allows great
flexibility in the choice of the threshold
parameter. A data-driven choice of the order threshold parameter is
given. An extension to a one-way HANOVA
setting is presented. Simulation studies with normal data suggest that
order thresholding can have great power
advantage over hard thresholding. Additional simulations with data
generated under a one-way HANOVA design
suggest even larger power gains over the traditional ANOVA $F$-test.

Applications of the order thresholding approach to testing for the
uniform distribution, and to multiple
testing problems will be pursued in a follow-up paper.

\begin{appendix}\label{app}
\section{\texorpdfstring{Proof of Theorem
\protect\lowercase{\ref{sing.seq.3.1}}}{Proof of Theorem 3.1}}

The proofs of the present lemmas can be found in
the archived supplemental material in \citet{KA09}.

\subsection{Some auxiliary results}

\begin{LEM}\label{aux.lem1}
Let $U_{i,n}$, $i=1,\ldots,n$, be order statistics from the uniform
distribution in $(0,1)$, and set
$V_{i,n}=-\log(1-U_{i,n})$. For any $0<\varepsilon<1$ and some
$1-\log(n-\sqrt{\frac{n}{2}\log(\frac{58}{\varepsilon
} )}+1 ) /\log
n\leq\delta(n)<1-\log(\frac{n}{2}\log(\frac
{58}{\varepsilon} ) ) /(2\log n)$, set
\[
\cases{
u_{jn}(\varepsilon)=\cases{
\displaystyle\max\Biggl\{0, \frac{j}{n}-\sqrt{\frac{1}{2n}
\log\biggl(\frac{58}{\varepsilon} \biggr)} \Biggr\},
&\quad $1\leq j<n^{1-\delta(n)}$,\vspace*{2pt}\cr
\displaystyle 1-e^{-\tilde{\nu}_{jn}}e^{\sqrt{2/\varepsilon}},
&\quad
$n^{1-\delta(n)}\leq j\leq n$,\cr}
\vspace*{2pt}\cr
u^{jn}(\varepsilon)=\cases{\displaystyle
\frac{j-1}{n}+\sqrt{\frac{1}{2n}\log\biggl(\frac{58}{\varepsilon
} \biggr)},
&\quad $1\leq j<n^{1-\delta(n)}$,\vspace*{2pt}\cr
\displaystyle 1-e^{-\tilde{\nu}_{jn}}e^{-\sqrt{2/\varepsilon}},
&\quad $n^{1-\delta(n)}\leq j\leq n$,}
}
\]
where $\tilde{\nu}_{jn}=\sum_{i=1}^j1/(n-i+1)$. Then, the sequences
of constants
\[
v_{jn}(\varepsilon)=-\log\bigl(1-u_{jn}(\varepsilon) \bigr),\qquad
v^{jn}(\varepsilon)=-
\log\bigl(1-u^{jn}(\varepsilon) \bigr)
\]
satisfy
%
%
\begin{equation} \label{def.v_jn}
P \{v_{jn}(\varepsilon)<V_{j,n}<v^{jn}(\varepsilon), 1\leq j\leq
n \}\geq1-\varepsilon,\qquad n\geq1.
\end{equation}
\end{LEM}
\begin{LEM}\label{aux.lem2}
Let $u_{jn}(\varepsilon)$ and $u^{jn}(\varepsilon)$ be given in Lemma
\ref
{aux.lem1}. Then, the sequences of
constants $u_{jn}(\varepsilon)$ and $u^{jn}(\varepsilon)$,
$j=1,\ldots,n$,
satisfy the relation
\[
u_{jn}(\varepsilon)<\frac{j}{n+1}<u^{jn}(\varepsilon).
\]
\end{LEM}
\begin{RE*} Assume that $1-n^{-\delta(n)}\to0$, as $n\to\infty$.
Then, the sequences of constants $u_{jn}(\varepsilon)$ and
$u^{jn}(\varepsilon)$, given in Lemma \ref{aux.lem1},
satisfy the relation $\sup_{1\leq j\leq n}(u^{jn}(\varepsilon
)-u_{jn}(\varepsilon))=o(1)$ (cf. Glivenko--Cantelli
theorem).
\end{RE*}
\begin{RE*} If we take all $u_{jn}(\varepsilon)$ and
$u^{jn}(\varepsilon)$
from the Kolmogorov's inequality, then $u_{jn}(\varepsilon
)=1-e^{-\tilde
\nu_{jn}+\sqrt{2/\varepsilon}}$ and
$u^{jn}(\varepsilon)=1-e^{-\tilde\nu_{jn}-\sqrt{2/\varepsilon}}$,
$j=1,\ldots,n$. Under these settings,
$\sup_{1\leq j\leq n}(u^{jn}(\varepsilon)-u_{jn}(\varepsilon))\neq o(1)$.
Also, the positive function~$R(j)$, defined
in Lemma \ref{aux.lem}, is not increasing on $1\leq j\leq n$.
\end{RE*}
\begin{LEM}\label{relation_vs} Let $\tilde{\nu}_{jn}$ be given in
Lemma \ref{aux.lem1}, and let
$\nu_{jn}=-\log(1-j/(n+1) )$, $j=1,\ldots,n$. Assume
that $1-n^{-\delta(n)}\to0$ and
$n^{1/2}(1-n^{-\delta(n)})\to\infty$, as $n\to\infty$. Then, the
sequences of constants $v_{jn}(\varepsilon)$ and
$v^{jn}(\varepsilon)$, given in Lemma \ref{aux.lem1}, satisfy the relations
\[
v_{jn}(\varepsilon)<\nu_{jn}<\tilde{\nu}_{jn}<v^{jn}(\varepsilon)
\quad\mbox{and}\quad
v^{jn}(\varepsilon)-v_{jn}(\varepsilon)\leq K(\varepsilon),
\]
where $K(\varepsilon)$ is independent of $n$.
\end{LEM}
\begin{LEM}\label{aux.lem} Let $v_{yn}(\varepsilon)$ and
$v^{yn}(\varepsilon)$, $1\leq y\leq
n$, be given in Lemma \ref{aux.lem1}, and let $\tilde H=F^{-1}\circ
G$, where $F$ is the central chi-squared
distribution function with 1 degree of freedom and $G$ is the standard
exponential distribution function. Then,
\begin{enumerate}
\item$\tilde H'$ is increasing, positive, concave, and $\tilde
H'(v)\to2$, as $v\to\infty$.
\item$\tilde H''$ is a decreasing positive function, and $\tilde
H''(v)\to0$, as $v\to\infty$.
\item$\tilde H(v)\tilde H''(v)\to0$, as $v\to\infty$.
\item$\frac{\tilde{H}'''(v)}{\tilde{H}''(v)} (1-e^{-v}
)\to
0$, as $v\to0$, and $\frac{\tilde{H}'''(v)}{\tilde{H}''(v)}\to
0$, as $v\to\infty$.
\item Assume
that $1-n^{-\delta(n)}\to0$ and $n^{1/2}(1-n^{-\delta(n)})\to\infty
$, as $n\to\infty$. The positive function
\[
R(y)=\bigl(v^{yn}(\varepsilon)-v_{yn}(\varepsilon)\bigr)\tilde{H}'' (
v_{yn}(\varepsilon) )\sqrt{\frac{y}{n-y+1}}
\]
is increasing on $1\leq y<n^{1-\delta(n)}$. Moreover, for sufficiently
large $n$, $R(y)$ is also increasing on
$n^{1-\delta(n)}\leq y\leq n$.
\end{enumerate}
\end{LEM}

\subsection{\texorpdfstring{Proof of Theorem
\protect\ref{sing.seq.3.1}}{Proof of Theorem 3.1}}

We need to check Assumptions A, B and C
of CGJ1967. We use the original forms of Assumptions A and C (restate
below for convenience), but a slightly stronger version of Assumption
B. [Note that the simultaneous bounds of
the exponential order statistics, $v_{jn}(\varepsilon)$ and
$v^{jn}(\varepsilon)$ used in Assumption B, are different
from those in CGJ1967.]

\textit{Assumption A}:
$\tilde H(v)$ is continuously differentiable for
$0<v<\infty$.

\textit{Assumption B}:
For each $\varepsilon>0$,
\begin{eqnarray*}
A_n&=&\sum_{j=n-k_n+1}^n \Biggl[ \Bigl\{{\sup_{v_{jn}(\varepsilon
)<v<v^{jn}(\varepsilon)}}
|\tilde{H}'(v)-\tilde{H}' (\tilde{\nu}_{jn}
) | \Bigr\} \sqrt{\frac{j}{n-j+1}} \Biggr]\\
&=&o(n\sigma_n(k_n)),
\end{eqnarray*}
where $v_{jn}(\varepsilon)$, $v^{jn}(\varepsilon)$, and $\tilde{\nu
}_{jn}$ are given in Lemma \ref{aux.lem1}.

\textit{Assumption C}:
${\max_{1\leq j\leq
n}} |\alpha_{jn}(k_n) |=o(n^{1/2}\sigma_n(k_n))$.

Assumption A is clearly satisfied. To verify Assumption C,
use Lemma \ref{aux.lem}(1) to write
%
%
\begin{eqnarray}\label{sing.seq.3.1.eq1}
\frac{{\max_{1\leq j\leq n}}|\alpha_{jn}(k_n)|}{\sqrt{n}\sigma
_n(k_n)}&=&\frac{\tilde{H}'
(\tilde{\nu}_{nn} )} {\sqrt{\sum_{j=1}^n\alpha
_{jn}^2(k_n)}}\nonumber\\[-8pt]\\[-8pt]
&\leq&
\frac{2}{\sqrt{\sum_{j=n-k_n+1}^n \{ \tilde{H}' (\tilde
{\nu}_{jn} ) \}^2}}.\nonumber
\end{eqnarray}
Suppose first that $k_n/n\rightarrow0$ as $n\rightarrow\infty$. Then
\[
\sum_{j=n-k_n+1}^n \{ \tilde{H}' (\tilde{\nu}_{jn}
) \}^2\geq
k_n \{\tilde{H}' (\tilde{\nu}_{n-k_n+1,n} ) \}
^2\to\infty\qquad \mbox{as } n\to\infty,
\]
so that (\ref{sing.seq.3.1.eq1}) tends to zero and Assumption C is
satisfied in this case. Next, suppose that
$k_n/n\to r$ as $n\rightarrow\infty$, for some $0<r\leq1$. Using the
approximation (2.9) of CGJ1967, that is,
$\tilde\nu_{jn}\simeq\nu_{jn}=-\log(1-j/(n+1) )$, it
follows that
\begin{eqnarray*}
&&\frac{1}{n}\sum_{j=n-k_n+1}^n \{ \tilde{H}' (\tilde{\nu
}_{jn} ) \}^2
\\
&&\qquad\simeq\frac{1}{n}\sum_{j=1}^n I \biggl(\frac{j}{n+1}>\frac
{n-k_n}{n+1} \biggr) \biggl\{
\biggl(1-\frac{j}{n+1} \biggr) (F^{-1} )' \biggl(\frac
{j}{n+1} \biggr) \biggr\}^2\\
&&\qquad\to\int_0^1 I(t>1-r) \biggl\{\frac{(1-t)}{f(F^{-1}(t))} \biggr\}
^2\,dt> 0\qquad \mbox{as }
n\to\infty,
\end{eqnarray*}
so that Assumption C is also satisfied. To show Assumption B, we use
Lemmas \ref{aux.lem}(1) and
\ref{relation_vs} to write ${\sup_{v_{jn}(\varepsilon
)<v<v^{jn}(\varepsilon)}}
|\tilde{H}'(v)-\tilde{H}' (\tilde{\nu}_{jn}
) | \leq
\tilde{H}'(v^{jn}(\varepsilon))-\tilde{H}' (v_{jn}(\varepsilon
) )=(v^{jn}
(\varepsilon)-v_{jn}(\varepsilon))\tilde{H}'' (\tilde
v_{jn}(\varepsilon
) )$, $\tilde v_{jn}(\varepsilon)\in
(v_{jn}(\varepsilon),v^{jn}(\varepsilon))$. Thus,
%
%
\begin{eqnarray}\label{num.1}\qquad\quad
\frac{A_n}{n\sigma_n(k_n)}&\leq&\frac{1}{\sqrt{n}}\sum
_{j=n-k_n+1}^n \Biggl[ \bigl\{
\bigl(v^{jn}(\varepsilon)-v_{jn}(\varepsilon)\bigr)\tilde{H}'' (
v_{jn}(\varepsilon) ) \bigr\}
\sqrt{\frac{j}{n-j+1}} \Biggr]\nonumber\\[-8pt]\\[-8pt]
&&{}\times\Biggl(\sqrt{\sum_{j=n-k_n+1}^n \{
\tilde{H}' (\tilde{\nu}_{jn} ) \}^2}\Biggr)^{-1},\nonumber
\end{eqnarray}
where the inequality is justified by the fact that $\tilde{H}''$ is a
decreasing positive function [Lemma
\ref{aux.lem}(2)]. We need to prove that (\ref{num.1}) tends to zero
as $n\to\infty$. Suppose first that
$k_n/n\to0$, as $n\rightarrow\infty$. Divide numerator and
denominator of (\ref{num.1}) by $k_n^{1/2}$ and
consider first the numerator. Then,
\begin{eqnarray*}
&&\hspace*{-6pt}\frac{1}{\sqrt{nk_n}}\sum_{j=n-k_n+1}^n \Biggl[ \bigl\{
\bigl(v^{jn}(\varepsilon)-v_{jn}(\varepsilon)\bigr)\tilde{H}'' (
v_{jn}(\varepsilon) ) \bigr\} \sqrt{\frac{j}{n-j+1}} \Biggr]\\
&&\hspace*{-6pt}\qquad=\cases{
\displaystyle\sqrt{\frac{8}{\varepsilon}}\dfrac{1}{\sqrt{nk_n}}\sum
_{j=n-k_n+1}^n \Biggl[
\tilde{H}'' (
v_{jn}(\varepsilon) )\sqrt{\frac{j}{n-j+1}} \Biggr],\cr
\qquad \mbox{if $\displaystyle k_n<\sqrt{\frac{n}{2}\log\biggl(\frac{58}{\varepsilon}
\biggr)}$},\cr
\displaystyle\frac{1}{\sqrt{nk_n}}\sum_{j=n-k_n+1}^{n^{1-\delta(n)}-1}
\Biggl[ \bigl\{
\bigl(v^{jn}(\varepsilon)-v_{jn}(\varepsilon)\bigr)\tilde{H}'' (
v_{jn}(\varepsilon) ) \bigr\}
\sqrt{\frac{j}{n-j+1}} \Biggr] \cr
\displaystyle\qquad{}+\frac{1}{\sqrt{nk_n}}\sum_{j=n^{1-\delta
(n)}}^{n-k_n^{1/4}} \Biggl[ \bigl\{
\bigl(v^{jn}(\varepsilon)-v_{jn}(\varepsilon)\bigr)\tilde{H}'' (
v_{jn}(\varepsilon) ) \bigr\} \sqrt{\frac{j}{n-j+1}} \Biggr]
\cr
\displaystyle\qquad{} +\frac{1}{\sqrt{nk_n}}\sum_{j=n-k_n^{1/4}+1}^n \Biggl[ \bigl\{
\bigl(v^{jn}(\varepsilon)-v_{jn}(\varepsilon)\bigr)\tilde{H}'' (
v_{jn}(\varepsilon) ) \bigr\} \sqrt{\frac{j}{n-j+1}} \Biggr],\cr
\qquad\mbox{otherwise},}
\end{eqnarray*}
with $1-\log(n-\sqrt{\frac{n}{2}\log(\frac{58}{\varepsilon
} )}+1 ) /\log
n\leq\delta(n)<1-\log(n-c_{\varepsilon}n^{3/16}k_n^{5/8}+1
) /\log n$ [This range is applied only
when $k_n\geq\sqrt{\frac{n}{2}\log(\frac{58}{\varepsilon
} )}$]. Assume that $1-n^{-\delta(n)}\to
0$, $n^{1/2}(1-n^{-\delta(n)})\to\infty$, and $n^{1/2}(1-n^{-\delta
(n)})^{3/2}\to d$ for some $d>0$, as
$n\to\infty$. If $k_n<\sqrt{\frac{n}{2}\log(\frac
{58}{\varepsilon} )}$, then
$\frac{1}{\sqrt{nk_n}}\sum_{j=n-k_n+1}^n [ \tilde{H}'' (
v_{jn}(\varepsilon) )\sqrt{\frac{j}{n-j+1}} ]
< (\frac{1}{2}\log(\frac{58}{\varepsilon} )
)^{1/4}n^{1/4} \tilde{H}'' (
v_{nn}(\varepsilon) ) \to0$, as $n\to\infty$. This inequality is
justified by Lemma \ref{aux.lem}(5), and
the fact that $n^{1/4} \tilde{H}'' ( v_{nn}(\varepsilon) )$
tends to
zero. 
Suppose that $k_n\geq\sqrt{\frac{n}{2}\log(\frac
{58}{\varepsilon} )}$ and $k_n/n\to0$, as $n\to\infty$.
Set\vspace*{2pt} $a_n=n^{1-\delta(n)}-1$ and $b_n=n-k_n^{1/4}$. Using Lemmas \ref
{aux.lem}(5) and \ref{relation_vs}, we
have
%
%
\begin{eqnarray}\label{pf1.1}
&&\frac{1}{\sqrt{nk_n}}\sum_{j=n-k_n+1}^n \Biggl[ \bigl\{
\bigl(v^{jn}(\varepsilon)-v_{jn}(\varepsilon)\bigr)\tilde{H}'' (
v_{jn}(\varepsilon) ) \bigr\} \sqrt{\frac{j}{n-j+1}}
\Biggr]\nonumber\\
&&\qquad\leq\sqrt{\frac{k_n}{n(1-n^{-\delta(n)})}} \bigl(v^{a_n,n}(\varepsilon
)-v_{a_n,n}(\varepsilon)\bigr)\tilde{H}'' (
v_{a_n,n}(\varepsilon) )\\
\label{pf1.1.1}
&&\qquad\quad{} + c_{\varepsilon}\sqrt{\frac{8}{\varepsilon
}}n^{3/16}\tilde{H}'' (
v_{b_n,n}(\varepsilon) )+\sqrt{\frac{8}{\varepsilon
}}k_n^{-1/4}\tilde{H}''(v_{nn}(\varepsilon)).
\end{eqnarray}
Since $(n-a_n+1)/n^{1/2}\to\infty$ as $n\to\infty$, using a
one-term Taylor expansion we have
$v^{a_n,n}(\varepsilon)-v_{a_n,n}(\varepsilon)\approx\frac{\sqrt
{2n\log
({58}/{\varepsilon} )}-1}
{n-a_n+1-\sqrt{{n}/{2}\log({58}/{\varepsilon}
)}}=O (\frac{1}{n^{1/2}(1-n^{-\delta(n)})} )$.
Thus, we have
\[
\mbox{(\ref{pf1.1})}=O \biggl(\frac{1}{n^{1/2}(1-n^{-\delta(n)})^{3/2}}\cdot
\frac{k_n^{1/2}}{n^{1/2}}\tilde{H}'' (
v_{a_n,n}(\varepsilon) ) \biggr)\to0\qquad \mbox{as } n\to\infty,
\]
where it is justified by Lemma \ref{aux.lem}(2) and the fact that
$v_{a_n,n}(\varepsilon)$ tends to infinity. From
Lemma \ref{aux.lem}(2), the second term of (\ref{pf1.1.1}) tends to 0 as
$n\to\infty$. Moreover, 
the first term of (\ref{pf1.1.1}) tends to 0 as $n\to\infty$ (even
$b_n=n-n^{1/4}$). Since also
$ (\frac{1}{k_n}\sum_{j=n-k_n+1}^n \{
\tilde{H}' (\tilde{\nu}_{jn} ) \}^2
)^{-1/2}\leq
(\tilde{H}'(\tilde{\nu}_{n-k_n+1,n}) )^{-1}<\infty$,
(\ref{num.1}) tends to zero and Assumption B is
satisfied when $k_n/n\to0$ as $n\to\infty$. Next, we suppose that
for some $0<r\leq1$, $k_n/n\to r$ as
$n\rightarrow\infty$. Divide numerator and denominator of (\ref
{num.1}) by $n^{1/2}$ and consider the numerator
and denominator separately. Since
\begin{eqnarray*}
&&\frac{1}{n}\sum_{j=n-k_n+1}^n \Biggl[ \bigl\{ \bigl(v^{jn}(\varepsilon
)-v_{jn}(\varepsilon)\bigr)\tilde{H}'' (
v_{jn}(\varepsilon) ) \bigr\} \sqrt{\frac{j}{n-j+1}} \Biggr]\\
&&\qquad\leq O \biggl(\frac{\tilde{H}'' (
v_{n^{1-\delta(n)}-1,n}(\varepsilon) )}{n^{1/2}(1-n^{-\delta
(n)})^{3/2}} \biggr)
+\sqrt{\frac{8}{\varepsilon}}n^{3/16}\tilde{H}'' (
v_{n-n^{1/4},n}(\varepsilon) )\\
&&\qquad\quad{}+\sqrt{\frac{8}{\varepsilon
}}n^{-1/4}\tilde{H}''(v_{nn}(\varepsilon))\\
&&\qquad\to0\qquad \mbox{as } n\to\infty,
\end{eqnarray*}
which can be obtained by breaking up the summation first for
$j=n-k_n+1$ to $n^{1-\delta(n)}-1$,
$n^{1-\delta(n)}$ to $n-n^{1/4}$, and lastly $n-n^{1/4}+1$ to $n$ with
$1-\log(n-\sqrt{\frac{n}{2}\log(\frac{58}{\varepsilon
} )}+1 ) /\log
n\leq\delta(n)<1-\log(n-n^{13/16}+1 ) /\log n$, and
$ (\frac{1}{n}\sum_{j=n-k_n+1}^n \{
\tilde{H}' (\tilde{\nu}_{jn} ) \}^2
)^{-1/2}<\infty$, the term (\ref{num.1}) converges to 0
as $n\to\infty$ in this case. Thus, Assumption B holds for both
cases. Since Assumptions A, B and C of CGJ1967
are satisfied, the proof is done.


\section{\texorpdfstring{Proof of Theorem
\protect\lowercase{\ref{anova.thm.4.1}}}{Proof of Theorem 4.1}}

The proofs of the present lemmas can be found in
the archived supplemental material in \citet{KA09}.

\subsection{Some auxiliary results}

\begin{LEM}\label{anova.u_v}
For any $0<\varepsilon<1$ and some $\delta(a)$ which satisfies
$1-\log(a-\sqrt{\frac{a}{2}\log(\frac{58}{\varepsilon
} )}+1 ) /\log
a\leq\delta(a)<1-\log(\frac{a}{2}\log(\frac
{58}{\varepsilon} ) ) /(2\log a)$,
$1-a^{-\delta(a)}\to0$, $a^{1/2}(1-a^{-\delta(a)})\to\infty$, and
$a^{1/2}(1-a^{-\delta(a)})^{3/2}\to d$ for
some $d>0$, as $a\to\infty$, let
\begin{eqnarray*}\label{const.u}
\cases{
u_{ja}(\varepsilon)=\cases{
\displaystyle\max\Biggl\{0, \frac{j}{a}-\sqrt{\frac{1}{2a}
\log\biggl(\frac{58}{\varepsilon} \biggr)} \Biggr\},
&\quad $1\leq j<a^{1-\delta(a)}$,\vspace*{2pt}\cr
\displaystyle1-e^{-\tilde{\nu}_{ja}}e^{\sqrt{2/\varepsilon}},
&\quad $a^{1-\delta(a)}\leq j\leq a$,}
\vspace*{2pt}\cr
u^{ja}(\varepsilon)=\cases{
\displaystyle\frac{j-1}{a}+\sqrt{\frac{1}{2a}\log\biggl(\frac{58}{\varepsilon
} \biggr)}, &\quad
$1\leq j<a^{1-\delta(a)}$,\vspace*{2pt}\cr
\displaystyle1-e^{-\tilde{\nu}_{ja}}e^{-\sqrt{2/\varepsilon}},
&\quad $a^{1-\delta(a)}\leq j\leq a$,}}
\end{eqnarray*}
where $\tilde{\nu}_{ja}=\sum_{i=1}^j1/(a-i+1)$, and set
\[
v_{ja}(\varepsilon)=-\log\bigl(1-u_{ja}(\varepsilon) \bigr),\qquad
v^{ja}(\varepsilon)=- \log\bigl(1-u^{ja}(\varepsilon) \bigr).
\]
For any $M\geq0$, let $\tilde H_{a,M}(v)=G_{a,M}^{-1}(1-e^{-v})$,
where $G_{a,M}$ is the noncentral
chi-squared distribution function with 1 degree of freedom and
noncentrality parameter $M^2/a$. Then,
\begin{enumerate}
%
\item$\tilde H_{a,M}'$ is bounded and $\tilde H_{a,M}'(v)\to2$,
as $v\to\infty$ and $a\to\infty$.
\item
$\tilde{H}_{a,M}''(v)=B_{a,M}(v)-(\tilde{H}_{a,M}'(v))^2J_{a,M}(v)$,
where
\begin{eqnarray*}
B_{a,M}(v)&=&-\tilde{H}_{a,M}'(v)+ (\tilde{H}_{a,M}'(v) )^2
\biggl(\frac{1}{2}+\frac{1}{2\tilde{H}_{a,M}(v)} \biggr),\\
J_{a,M}(v)&=&\frac{\sum_{k=1}^{\infty} \{{(M^2/a)^k
(\tilde{H}_{a,M}(v) )^{k-1}}/
({2^{2k}(k-1)!\Gamma(k+1/2 )}) \}}
{\sum_{k=0}^{\infty} \{{(M^2/a)^k (\tilde
{H}_{a,M}(v) )^k}/
({2^{2k}k!\Gamma(k+1/2 )}) \}}.
\end{eqnarray*}
Note that $B_{a,M}$ is a decreasing positive function, $B_{a,M}(v)\to
0$ as $v\to\infty$ and $a\to\infty$, and
$J_{a,M}$ is bounded by $M^2/(2a)$.
\item The positive function
\[
R_M(y)= \bigl(v^{ya}(\varepsilon)-v_{ya}(\varepsilon) \bigr)
B_{a,M}(v_{ya}(\varepsilon))\sqrt{\frac{y}{a-y+1}}
\]
is increasing on $1\leq y<a^{1-\delta(a)}$. Moreover, for sufficiently
large $a$, $R_M(y)$ is also increasing
on $a^{1-\delta(a)}\leq y\leq a$.
\end{enumerate}
\end{LEM}
\begin{LEM}\label{rel.chi-nonchi}
Consider the setting of Lemma \ref{anova.u_v}. Let $g_{a,0}$ and
$g_{a,M}$ be the density functions of
$\chi_1^2(0)$ and $\chi^2_1(M^2/a)$, respectively. Set
$y_{a,0}=\tilde
H_{a,0}(v_{aa}(\varepsilon))=G_{a,0}^{-1}(1-e^{-v_{aa}(\varepsilon
)})$ and
$y_{a,M}=\tilde
H_{a,M}(v_{aa}(\varepsilon))=G_{a,M}^{-1}(1-e^{-v_{aa}(\varepsilon
)})$. Then,
\begin{enumerate}
\item$y_{a,M}$ is bounded by $ (2\log(a+1)-2\log(\sqrt
{\pi/2} )+2\log(
e^{M/(2\sqrt{a})}+\break e^{-M/(2\sqrt{a})} ) )^2$.
\item$g_{a,0}(y_{a,M})/g_{a,M}(y_{a,M})\to1$ and $a^{1/4}
(g_{a,0}(y_{a,M})/g_{a,M}(y_{a,M})-1 )\to
0$, as $a\to\infty$.
\item$\frac
{g_{a,0}(y_{a,0})-g_{a,0}(y_{a,M})}{g_{a,0}(y_{a,M})}\approx
- (\frac{y_{a,M}^{-1}+1}{2} )\frac{M^2}{2a}C_{a,M}$,
where $C_{a,M}$ is defined in the proof. In
particular, $C_{a,M}=O(y_{a,M})$.
\item$a^{1/4} (y_{a,0}/y_{a,M}-1 )\to
0$ and $a^{1/4} (g_{a,0}(y_{a,0})/g_{a,0}(y_{a,M})-1 )
\to0$, as $a\to\infty$.
\item$a^{1/4} (e^{-v_{aa}(\varepsilon
)}/g_{a,0}(y_{a,M})-2+2/y_{a,0} )\to
0$, as $a\to\infty$.
\item$a^{1/4} (e^{-v_{aa}(\varepsilon
)}/g_{a,M}(y_{a,M})-2+2/y_{a,M} )\to
0$, as $a\to\infty$.
\end{enumerate}
\end{LEM}
\begin{RE*} For any $M\geq0$, we write
$\tilde H'_{a,M}(v_{aa}(\varepsilon))=\frac{e^{-v_{aa}(\varepsilon)}}
{g_{a,M}(y_{a,M})}$ and\break $B_{a,M}(v_{aa}(\varepsilon))
=\frac{e^{-v_{aa}(\varepsilon)}}{2g_{a,M} (y_{a,M})}$ $ (
\frac{e^{-v_{aa}(\varepsilon)}}{g_{a,M}(y_{a,M})}-2+\frac
{e^{-v_{aa}(\varepsilon)}}{g_{a,M}(y_{a,M})}
\cdot\frac{1}{y_{a,M}} )$. From Lemmas~\ref{anova.u_v}(1),
\ref{anova.u_v}(2) and \ref{rel.chi-nonchi}(6), we obtain
that $a^{1/4}B_{a,M}(v_{aa}(\varepsilon))\to0$, as $a\to\infty$.
\end{RE*}
\begin{LEM}\label{rel.chi-nonchi.add}
Consider the setting of Lemma \ref{rel.chi-nonchi}. Let
$b_a=a-k_a^{1/4}$ with
$k_a\geq\sqrt{\frac{a}{2}\log(\frac{58}{\varepsilon}
)}$. Set $x_{a,0}=\tilde
H_{a,0}(v_{b_a,a}(\varepsilon
))=G_{a,0}^{-1}(1-e^{-v_{b_a,a}(\varepsilon
)})$ and $x_{a,M}=\tilde
H_{a,M}(v_{b_a,a}(\varepsilon
))=G_{a,M}^{-1}(1-e^{-v_{b_a,a}(\varepsilon
)})$. Then,
\begin{enumerate}
\item$x_{a,M}\leq(2\log(a+1)-2\log(k_a^{1/4}+1)-2\log
(\sqrt{\pi/2} )+2\log(
e^{M/(2\sqrt{a})}+e^{-M/(2\sqrt{a})} ) )^2$.
\item$g_{a,0}(x_{a,M})/g_{a,M}(x_{a,M})\to1$ and $a^{3/16}
(g_{a,0}(x_{a,M})/g_{a,M}(x_{a,M})-1 )\to
0$, as $a\to\infty$.
\item$\frac
{g_{a,0}(x_{a,0})-g_{a,0}(x_{a,M})}{g_{a,0}(x_{a,M})}\approx
- (\frac{x_{a,M}^{-1}+1}{2} )\frac{M^2}{2a}C'_{a,M}$, where
\[
C'_{a,M}=\frac{\phi' (\sqrt{x_{a,M}}-t^*/\sqrt{a} )
+\phi' (\sqrt{x_{a,M}}+t^*/\sqrt{a}
)}{g_{a,0}(x_{a,\tilde t})} \qquad\mbox{with } t^*,\tilde t\in
(0,M).
\]
In particular, $C'_{a,M}=O(x_{a,M})$.
\item$a^{3/16} (x_{a,0}/x_{a,M}-1 )\to
0$ and $a^{3/16} (g_{a,0}(x_{a,0})/g_{a,0}(x_{a,M})-1 )
\to0$, as $a\to\infty$.
\item$a^{3/16} (e^{-v_{b_a,a}(\varepsilon
)}/g_{a,0}(x_{a,M})-2+2/x_{a,0} )\to
0$, as $a\to\infty$.
\item$a^{3/16} (e^{-v_{b_a,a}(\varepsilon
)}/g_{a,M}(x_{a,M})-2+2/x_{a,M} )\to
0$, as $a\to\infty$.
\end{enumerate}
\end{LEM}
\begin{RE*} From Lemmas \ref{anova.u_v}(1), \ref{anova.u_v}(2),
and \ref{rel.chi-nonchi.add}(6), we obtain that
$a^{3/16}\times B_{a,M}(v_{a-k_a^{1/4},a}(\varepsilon))\to0$, as
$a\to\infty$.
\end{RE*}

\subsection{\texorpdfstring{Proof of Theorem
\protect\ref{anova.thm.4.1}}{Proof of Theorem 4.1}}

For simplicity, let $\tilde{H}_{a,t}(v)=G_{a,t}^{-1}(1-e^{-v})$. Then
\[
\alpha_{ia}^t(k_a)=\frac{1}{a-i+1}\sum_{j=i}^ac_{ja}\frac{
e^{-\tilde{\nu}_{ja}}}
{g_{a,t} (G_{a,t}^{-1}(1-e^{-\tilde{\nu}_{ja}}) )}
=
\frac{1}{a-i+1}\sum_{j=i}^ac_{ja}\tilde{H}_{a,t}' (\tilde{\nu
}_{ja} )
\]
and
\[
(\sigma_a^t(k_a) )^2=\frac{1}{a}\sum_{i=1}^a(\alpha
_{ia}^t(k_a))^2.
\]
Let us check Assumptions A, B and C of CGJ1967, which we restated in
the proof of Theorem \ref{sing.seq.3.1}.
For given any $|t|<M$, Assumption A is clearly satisfied. Next, it is
easily verified that for any fixed values
of $a$ and $v$, $\tilde H'_{a,t}(v)$ increases as $|t|$ increases.
Thus, $\alpha_{ia}^t(k_a)$ and
$\sigma_a^t(k_a)$ increase as $|t|$ increases. Let us check Assumption
C: for given any $|t|<M$,
\begin{eqnarray*}
\frac{{\max_{1\leq j\leq a}} |\alpha_{ja}^t(k_a) |}{\sqrt
{a}\sigma_a^t(k_a)}&\leq&\frac{{\max_{1\leq
j\leq a}} |\alpha_{ja}^M(k_a) |}{\sqrt{a}\sigma_a^0(k_a)}
\\
&\leq&\frac{\max_{a-k_a+1\leq j\leq
a}\tilde{H}'_{a,M} (\tilde{\nu}_{ja} )} {\sqrt{\sum
_{j=a-k_a+1}^a \{\tilde
H_{a,0}'(\tilde{\nu}_{ja}) \}^2}}\\
&\to&0 \qquad
\mbox{as } a\to \infty,
\end{eqnarray*}
provided that $k_a\to\infty$, as $a\to\infty$. It is justified by
the facts that
\[
\max_{a-k_a+1\leq j\leq
a}\tilde{H}'_{a,M} (\tilde{\nu}_{ja} )
\]
is bounded [Lemma \ref{anova.u_v}(1)] and
$\sum_{j=a-k_a+1}^a \{\tilde H_{a,0}'(\tilde{\nu}_{ja}) \}
^2\to\infty$ as $k_a$ tends to infinity with
$a$. (It was shown in the proof of Theorem \ref{sing.seq.3.1} because
it becomes the central chi-square case
when $t=0$.) In order to verify Assumption B, it suffices to show that
\[
\frac{\sum_{j=a-k_a+1}^a [ \{{\sup_{v_{ja}(\varepsilon
)<v<v^{ja}(\varepsilon)}}
|\tilde{H}_{a,M}'(v)-\tilde{H}_{a,M}' (\tilde{\nu
}_{ja} ) | \}
\sqrt{{j}/({a-j+1})} ]}{\sqrt{a\sum_{j=a-k_a+1}^a \{
\tilde{H}_{a,0}'
(\tilde{\nu}_{ja} ) \}^2}}=o(1),
\]
where $v_{ja}(\varepsilon)$, $v^{ja}(\varepsilon)$, and $\tilde\nu_{ja}$
are given in Lemma \ref{anova.u_v}. Using
Lemma \ref{anova.u_v}(2), we write
\begin{eqnarray*}
&&{\sup_{v_{ja}(\varepsilon)<v<v^{ja}(\varepsilon)}}
|\tilde{H}_{a,M}'(v)-\tilde{H}_{a,M}' (\tilde{\nu
}_{ja} ) |\\
&&\qquad\leq\bigl(v^{ja}(\varepsilon)-v_{ja}(\varepsilon) \bigr)\cdot
|\tilde
H_{a,M}''(v_{ja}^*) |\\
&&\qquad \leq\bigl(v^{ja}(\varepsilon)-v_{ja}(\varepsilon)
\bigr)B_{a,M}(v_{ja}(\varepsilon))\\
&&\qquad\quad{}+ \bigl(v^{ja}(\varepsilon)-v_{ja}(\varepsilon) \bigr) (\tilde
{H}_{a,M}'(v_{ja}^*) )^2M^2/(2a)
\end{eqnarray*}
with some $v_{ja}^*\in(v_{ja}(\varepsilon),v^{ja}(\varepsilon))$.
From the
above inequality, we have
%
%
\begin{eqnarray}\hspace*{5pt}
\label{subproof-2.9}
&&\frac{\sum_{j=a-k_a+1}^a [ \{{\sup_{v_{ja}(\varepsilon
)<v<v^{ja}(\varepsilon)}}
|\tilde{H}_{a,M}'(v)-\tilde{H}_{a,M}' (\tilde{\nu
}_{ja} ) | \}
\sqrt{{j}/({a-j+1})} ]}{\sqrt{a\sum_{j=a-k_a+1}^a \{
\tilde{H}_{a,0}'
(\tilde{\nu}_{ja} ) \}^2}}\nonumber\\
&&\qquad\leq\frac{{1}/{\sqrt{a}}\sum_{j=a-k_a+1}^a [ \{
(v^{ja}(\varepsilon)-v_{ja}(\varepsilon) )
B_{a,M}(v_{ja}(\varepsilon)) \} \sqrt{{j}/({a-j+1})}
]}{\sqrt{\sum_{j=a-k_a+1}^a \{\tilde{H}_{a,0}'
(\tilde{\nu}_{ja} ) \}^2}}\\
\label{subproof-2.9.1}
&&\qquad\quad{}
+\frac{1}{\sqrt{a}}
\sum_{j=a-k_a+1}^a \Biggl[ \biggl\{ \bigl(v^{ja}(\varepsilon
)-v_{ja}(\varepsilon) \bigr)
(\tilde{H}_{a,M}'(v_{ja}^*) )^2\frac{M^2}{2a} \biggr\}
\sqrt{\frac{j}{a-j+1}} \Biggr]\nonumber\\[-8pt]\\[-8pt]
&&\qquad\quad\hspace*{11pt}{}\times\Biggl({\sqrt{\sum_{j=a-k_a+1}^a \{
\tilde{H}_{a,0}'
(\tilde{\nu}_{ja} ) \}^2}}\Biggr)^{-1}.\nonumber
\end{eqnarray}
To show that (\ref{subproof-2.9.1}) tends to zero, we use Lemmas \ref
{anova.u_v}(1) and \ref{relation_vs}
to write
\[
\mbox{(\ref{subproof-2.9.1})}\leq C_{\varepsilon}\cdot\frac
{k_a}{a}\cdot
\Biggl(\sum_{j=a-k_a+1}^a \{\tilde{H}_{a,0}' (\tilde{\nu
}_{ja} ) \}^2 \Biggr)^{-1/2}
\qquad\mbox{for some } 0<C_{\varepsilon}<\infty.
\]
Suppose first that $k_a/a\to0$ as $a\to\infty$. Then
%
%
\begin{equation}\label{ris0}
\Biggl(\frac{1}{k_a}\sum_{j=a-k_a+1}^a \{\tilde{H}_{a,0}'
(\tilde{\nu}_{ja} ) \}^2 \Biggr)^{-1/2}\leq\frac
{1}{\tilde
H'_{a,0}(\tilde{\nu}_{a-k_a+1,a})}<\infty,
\end{equation}
so that (\ref{subproof-2.9.1}) tends to 0 as $a\to\infty$. For some
$0<r\leq1$, if $k_a/a\to r$ as
$a\to\infty$, then
%
%
\begin{equation}\label{ris1}
\Biggl(\frac{1}{a}\sum_{j=a-k_a+1}^a \{\tilde{H}_{a,0}'
(\tilde{\nu}_{ja} ) \}^2 \Biggr)^{-1/2}<\infty.
\end{equation}
Thus, (\ref{subproof-2.9.1}) tends to 0 as $a\to\infty$ in both
cases. Since (\ref{subproof-2.9.1}) converges to
zero, the remaining part is to prove that (\ref{subproof-2.9}) tends
to 0, provided that $k_a\to\infty$ as
$a\to\infty$. Suppose first that $k_a<\sqrt{\frac{a}{2}\log
(\frac{58}{\varepsilon} )}$. Divide numerator
and denominator of (\ref{subproof-2.9}) by $k_a^{1/2}$ and consider
first the numerator. From Lemma
\ref{anova.u_v}(3), we have
\begin{eqnarray*}
&&\frac{1}{\sqrt{ak_a}}\sum_{j=a-k_a+1}^a \Biggl[ \bigl\{
\bigl(v^{ja}(\varepsilon)-v_{ja}(\varepsilon)\bigr)B_{a,M}(v_{ja}(\varepsilon))
\bigr\} \sqrt{\frac{j}{a-j+1}} \Biggr]\\
&&\qquad\leq\sqrt{\frac{8}{\varepsilon}} \biggl(\frac{1}{2}\log
\biggl(\frac{58}{\varepsilon} \biggr) \biggr)^{1/4}
a^{1/4}B_{a,M} (v_{aa}(\varepsilon) ).
\end{eqnarray*}
Using Lemmas \ref{anova.u_v}(1), \ref{anova.u_v}(2), \ref
{rel.chi-nonchi}(6) and (\ref{ris0}), the term
(\ref{subproof-2.9}) tends to zero and Assumption B is satisfied in
this case. Next, we suppose that
$k_a\geq\sqrt{\frac{a}{2}\log(\frac{58}{\varepsilon} )}$
and $k_a/a\to0$, as $a\to\infty$. Then, from
Lemmas \ref{anova.u_v}(3) and \ref{relation_vs},
\begin{eqnarray*}
&&\frac{1}{\sqrt{ak_a}}\sum_{j=a-k_a+1}^a \Biggl[ \bigl\{
\bigl(v^{ja}(\varepsilon)-v_{ja}(\varepsilon)\bigr)B_{a,M}(v_{ja}(\varepsilon
)) \bigr\}
\sqrt{\frac{j}{a-j+1}} \Biggr]\\
&&\qquad\leq\sqrt{\frac{k_a}{a(1-a^{-\delta(a)})}}\bigl(v^{a^{1-\delta
(a)}-1,a}(\varepsilon)-v_{a^{1-\delta(a)}-1,a}(\varepsilon)\bigr)
B_{a,M} (v_{a^{1-\delta(a)}-1,a}(\varepsilon) )\\
&&\qquad\quad{}+c_{\varepsilon}\sqrt{\frac{8}{\varepsilon}}a^{3/16}B_{a,M} (
v_{a-k_a^{1/4},a}(\varepsilon) )+\sqrt{\frac{8}{\varepsilon
}}k_a^{-1/4}B_{a,M}(v_{aa}(\varepsilon)),
\end{eqnarray*}
where $1-\log(a-\sqrt{\frac{a}{2}\log(\frac
{58}{\varepsilon} )}+1 ) /\log
a\leq\delta(a)<1-\log(a-c_{\varepsilon}a^{3/16}k_a^{5/8}+1
) /\log a$. From Lemmas \ref{anova.u_v}(1),
\ref{anova.u_v}(2), \ref{rel.chi-nonchi.add}(6), (\ref{ris0}), and
the fact that\break
$v^{a^{1-\delta(a)}-1,a}(\varepsilon)-v_{a^{1-\delta
(a)}-1,a}(\varepsilon)
=O (\frac{1}{a^{1/2}(1-a^{-\delta(a)})} )$, the term (\ref
{subproof-2.9}) also tends to zero and
Assumption B is satisfied in this case. Lastly, we suppose that for
some $0<r\leq1$, $k_a/a\to r$ as
$a\to\infty$. Divide numerator and denominator of (\ref
{subproof-2.9}) by $a^{1/2}$ and consider the numerator
and denominator separately. Since (\ref{ris1}) and
\begin{eqnarray*}
&&\frac{1}{a}\sum_{j=a-k_a+1}^a \Biggl[ \bigl\{ \bigl(v^{ja}(\varepsilon
)-v_{ja}(\varepsilon)\bigr)B_{a,M} (
v_{ja}(\varepsilon) ) \bigr\} \sqrt{\frac{j}{a-j+1}} \Biggr]\\
&&\qquad\leq\frac{1}{\sqrt{1-a^{-\delta(a)}}}\bigl(v^{a^{1-\delta
(a)}-1,a}(\varepsilon)-v_{a^{1-\delta(a)}-1,a}(\varepsilon)\bigr)
B_{a,M} (v_{a^{1-\delta(a)}-1,a}(\varepsilon) )\\
&&\qquad\quad{}+\sqrt{\frac{8}{\varepsilon}}a^{3/16}B_{a,M} (
v_{a-a^{1/4},a}(\varepsilon) )+\sqrt{\frac{8}{\varepsilon
}}a^{-1/4}B_{a,M}(v_{aa}(\varepsilon))\\
&&\qquad\to0 \qquad\mbox{as } a\to\infty
\end{eqnarray*}
with $1-\log(a-\sqrt{\frac{a}{2}\log(\frac{58}{\varepsilon
} )}+1 ) /\log
a\leq\delta(a)<1-\log(a-a^{13/16}+1 ) /\log a$, the
term (\ref{subproof-2.9}) converges to 0 as
$a\to\infty$ in this case. Thus, Assumption B holds as $k_a$ tends to
infinity with $a$. Since Assumptions A,
B and C of CGJ1967 are satisfied, the proof is done.


\section{\texorpdfstring{Proof of Theorem
\protect\lowercase{\ref{anova.t.random}}}{Proof of Theorem 4.2}}
\subsection{\texorpdfstring{Proof of Lemmas
\protect\ref{anova.lem.4.2}--\protect\ref{anova.thm.4.2}}{Proof of Lemmas 4.2--4.4}}
\subsubsection{\texorpdfstring{Proof of Lemma
\protect\ref{anova.lem.4.2}}{Proof of Lemma 4.2}}

We observe that
\[
\frac{\sqrt{a}Q_a^t(k_a)}{\sigma_a^t(k_a)}=\frac{1}{\sqrt
{a}\sigma_a^t(k_a)}
\sum_{i=1}^a\alpha_{ia}^t(k_a)(V_i-1),
\]
where $V_i$ are i.i.d. random variables with the distribution function
$G(v)=1-e^{-v}$, $v\geq0$. Note that
$g(v)=e^{-v}$, $v\geq0$, $E(V_i-1)=0$, $\operatorname{Var}(V_i-1)=1$, and $ E
(|V_i-1|^3 )=12/e-2$. Let
\begin{eqnarray*}
s_a^2&=&\sum_{i=1}^a\operatorname{Var} \bigl(\alpha_{ia}^t(k_a)(V_i-1) \bigr)=\sum
_{i=1}^a (\alpha_{ia}^t(k_a) )^2
=a (\sigma_a^t(k_a) )^2,\\
\beta_a^3&=&\sum_{i=1}^aE \bigl(|\alpha_{ia}^t(k_a)(V_i-1)|^3
\bigr)= \biggl(\frac{12}{e}-2 \biggr)
\sum_{i=1}^a|\alpha_{ia}^t(k_a)|^3,\\
r_a&=&\frac{\beta_a^3}{s_a^3}= \biggl(\frac{12}{e}-2 \biggr)\frac
{\sum_{i=1}^a|\alpha_{ia}^t(k_a)|^3}
{a^{3/2} (\sigma_a^t(k_a) )^3}.
\end{eqnarray*}
Using Berry--Esseen theorem of Galambos [(\citeyear{G95}), page 180] we have
\[
{\sup_{-\infty< x<\infty}} |S_a^t (xs_a )-\Phi
(x) |\leq0.8r_a\qquad \mbox{as }
a\rightarrow\infty,
\]
where $S_a^t$ is a distribution function of $\sum_{i=1}^a\alpha
_{ia}^t(k_a)(V_i-1)$ and $\Phi$ is a standard
normal distribution function. Thus, we have
\begin{eqnarray*}
&&{\mathop{\sup_{-M< t< M}}_{-\infty< x<\infty}} |H_{a,t}(x)-\Phi
(x)|\\
&&\qquad=
{\mathop{\sup_{-M< t< M}}_{-\infty< x<\infty}} \bigl|S_a^t
\bigl(x\sqrt{a}\sigma_a^t(k_a) \bigr)-\Phi(x) \bigr|\\
&&\qquad\leq0.8 \biggl(\frac{12}{e}-2 \biggr)\sup_{-M< t< M} \biggl\{
\frac{{\max_{1\leq i\leq
a}}|\alpha_{ia}^t(k_a)|} {\sqrt{a}\sigma_a^t(k_a)} \biggr\}\to0
\qquad\mbox{as } a\to\infty,
\end{eqnarray*}
provided that $k_a\to\infty$, as $a\to\infty$.

\subsubsection{\texorpdfstring{Proof of Lemma
\protect\ref{anova.lem.4.3}}{Proof of Lemma 4.3}}

For convenience, we rewrite
\[
R_a^t(k_a)=\frac{1}{a}\sum_{j=a-k_a+1}^a \{
(V_{j,a}-\tilde{\nu}_{ja} )G_{ja}^t (V_{j,a}
) \},
\]
where
\[
G_{ja}^t(v)=\cases{
\dfrac{\tilde{H}_{a,t}(v)-\tilde{H}_{a,t}(\tilde{\nu
}_{ja})}{v-\tilde{\nu}_{ja}}
-\tilde{H}_{a,t}' (\tilde{\nu}_{ja} ), &\quad if
$v\neq\tilde{\nu}_{ja}$,\cr
0, &\quad if $v=\tilde{\nu}_{ja}$.}
\]
Let $g_{ja}^t(\varepsilon)={\sup_{v_{ja}(\varepsilon
)<v<v^{ja}(\varepsilon
)} }|G_{ja}^t(v) |$, where
$v_{ja}(\varepsilon)$ and $v^{ja}(\varepsilon)$ are given in Lemma
\ref
{anova.u_v}. Then, we have
%
%
\begin{equation}\label{20080917-1.1}\qquad
P \Biggl\{ |R_a^t(k_a) |\leq\frac{1}{a}\sum_{j=a-k_a+1}^a
g_{ja}^M(\varepsilon) |V_{j,a}-\tilde{\nu}_{ja} | \mbox
{ for all } |t|<M \Biggr\}\geq1-\varepsilon.
\end{equation}
It follows from (\ref{20080917-1.1}) that
\[
P \Biggl\{{\sup_{-M<t<M} }|R_a^t(k_a) |\leq\frac{1}{a}\sum
_{j=a-k_a+1}^a
g_{ja}^M(\varepsilon) |V_{j,a}-\tilde{\nu}_{ja} | \Biggr\}
\geq1-\varepsilon.
\]
From Assumption B and Proposition 2 of CGJ1967, we have\break $\sum_{j=a-k_a+1}^a
g_{ja}^M(\varepsilon) \times|V_{j,a}-\tilde{\nu}_{ja} |=o_p
(\sqrt{a}\sigma_a^M(k_a) )$, so that\break
$\sup_{-M<t<M} |\sqrt{a}R_a^t(k_a) |=o_p (\sigma
_a^M(k_a) )$. Also, it is easily verified
that\break $\sigma_a^M(k_a)/\sigma_a^0(k_a)=O(1)$ (Lemma \ref{ratio.std}),
provided that $k_a\to\infty$ as
$a\to\infty$. Consequently,
\[
\sup_{-M<t<M} \biggl|\frac{\sqrt{a}R_a^t(k_a)}{\sigma
_a^t(k_a)} \biggr|\leq
\sup_{-M<t<M}\frac{ |\sqrt{a}R_a^t(k_a) |}{\sigma
_a^0(k_a)}=o_p(1).
\]

\subsubsection{\texorpdfstring{Proof of Lemma
\protect\ref{anova.thm.4.2}}{Proof of Lemma 4.4}}

We have already proved that for given any $|t|<M$,
\begin{eqnarray*}
{T_L^{t}}^*(k_a)&=&\frac{\sqrt{a}Q_a^t(k_a)}{\sigma_a^t(k_a)}\\
&&{}+\frac
{\sqrt{a}R_a^t(k_a)}{\sigma_a^t(k_a)}
\stackrel{d}{\rightarrow}N(0,1)\qquad \mbox{as } a\to\infty\qquad
(\mbox{Theorem \ref{anova.thm.4.1}})
\end{eqnarray*}
provided that $k_a\to\infty$, as $a\to\infty$. From Lemmas \ref
{anova.lem.4.2} and \ref{anova.lem.4.3}, we have
\[
{\mathop{\sup_{-M< t< M}}_{-\infty< x<\infty}}
|F_{a,t}(x)-\Phi(x) |\rightarrow0\qquad \mbox{as }
a\rightarrow\infty,
\]
provided that $k_a\to\infty$, as $a\to\infty$.

\subsection{\texorpdfstring{Proof of Theorem
\protect\ref{anova.t.random}}{Proof of Theorem 4.2}}

For any given $\delta_1>0$, there exists $M>0$ such that
\[
P ( |t |\geq M )<\delta_1.
\]
From Lemma \ref{anova.thm.4.2}, any given $\delta_2>0$, there exists
$a_0$ such that
\[
\bigl|P \bigl({\widehat T_L}^*(k_a)\leq x | |t |<
M \bigr)-\Phi(x) \bigr|<\delta_2\qquad\mbox{for
all }a>a_0.
\]
Thus, we have
\begin{eqnarray*}
\bigl|P \bigl({\widehat T_L}^*(k_a)\leq x \bigr)-\Phi(x)\bigr|
&\leq&\bigl|P \bigl({\widehat T_L}^*(k_a)\leq
x | |t |< M \bigr)-\Phi(x) \bigr|\\
&&{}+P (
|t |\geq M ) <\delta_1+\delta_2
\end{eqnarray*}
for all $a>a_0$. Take $\varepsilon/2=\max\{\delta_1,\delta
_2 \}$. Then
\[
\bigl|P \bigl({\widehat T_L}^*(k_a)\leq x \bigr)-\Phi(x)
\bigr|<\varepsilon\qquad\mbox{for all } a>a_0.
\]
Thus, provided that $k_a\to\infty$, as $a\to\infty$, we have
\[
{\widehat T_L}^*(k_a)\stackrel{d}{\rightarrow}N(0,1)\qquad
\mbox{as } a\to\infty.
\]


\section{\texorpdfstring{Proof of Theorem
\protect\lowercase{\ref{anova.order.thres}}}{Proof of Theorem 4.3}}
\subsection{\texorpdfstring{Proof of Lemmas \protect\ref{ratio.std}--\protect\ref
{second.lem}}{Proof of Lemmas 4.5--4.7}}

\subsubsection{\texorpdfstring{Proof of Lemma \protect\ref{ratio.std}}{Proof of Lemma 4.5}} We need to show
that
\[
\sup_{-M< t< M}
\biggl|\frac{\sigma_a^t(k_a)}{\sigma_a^0(k_a)}-1 \biggr|=
\frac{\sigma_a^M(k_a)-\sigma_a^0(k_a)}{\sigma_a^0(k_a)}
\to0\qquad \mbox{as } a\to\infty,
\]
provided that $k_a\to\infty$ as $a\to\infty$. Suppose first that
$k_a/a\to0$ as $a\to\infty$. Since $\sqrt{a/k_a}\sigma_a^0(k_a)>0$,
it is enough to show that $
\sqrt{a/k_a} (\sigma_a^M(k_a)-\sigma_a^0(k_a) )\to0$, as
$a\to\infty$. We first have
$a(\sigma_a^M(k_a))^2\leq(k_a (2-\break k_a/a ) )
\{\max_{a-k_a+1\leq j\leq a}\tilde
H_{a,M}'(\tilde{\nu}_{ja}) \}^2$ and $a(\sigma_a^0(k_a))^2\geq
(k_a [1+k_a(a-k_a)/
((a+1)(k_a+1)) ] ) \{\tilde
H_{a,0}'(\tilde{\nu}_{a-k_a+1,a}) \}^2$. Consequently, we
obtain
\begin{eqnarray*}
&&\sqrt{\frac{a}{k_a}} \bigl(\sigma_a^M(k_a)-\sigma_a^0(k_a)
\bigr)\\
&&\qquad\leq
\sqrt{2-\frac{k_a}{a}} \Bigl\{\max_{a-k_a+1\leq j\leq a}\tilde
H_{a,M}'(\tilde{\nu}_{ja}) \Bigr\}\\
&&\qquad\quad{}-\sqrt{1+\frac
{k_a(a-k_a)}{(a+1)(k_a+1)}} \{\tilde
H_{a,0}'(\tilde{\nu}_{a-k_a+1,a}) \}\to0
\end{eqnarray*}
as $a\to\infty$. Next, we suppose that for some $0<r\leq1$,
$k_a/a\to r$ as $a\to\infty$. Then
$\sigma_a^0(k_a)>0$, so we need to prove that $\sigma_a^M(k_a)-\sigma
_a^0(k_a)\to0$, as $a\to\infty$. We
observe that
\begin{eqnarray*}
&&(\sigma_a^M(k_a))^2-(\sigma_a^0(k_a))^2\\
&&\qquad=\frac{1}{a}\sum_{i=1}^{a-k_a} \biggl(\frac{1}{a-i+1}
\biggr)^2 \Biggl\{ \Biggl(\sum_{j=a-k_a+1}^a \tilde
H_{a,M}'(\tilde{\nu}_{ja}) \Biggr)^2\\
&&\qquad\quad\hspace*{95.9pt}{} - \Biggl(\sum_{j=a-k_a+1}^a
\tilde H_{a,0}'(\tilde{\nu}_{ja}) \Biggr)^2 \Biggr\}\\
&&\qquad\quad{}+\frac{1}{a}\sum_{i=a-k_a+1}^a \Biggl\{ \Biggl(\frac
{1}{a-i+1}\sum_{j=i}^a \tilde
H_{a,M}'(\tilde{\nu}_{ja}) \Biggr)^2\\
&&\qquad\quad\hspace*{68.2pt}{} - \Biggl(\frac{1}{a-i+1}\sum
_{j=i}^a \tilde
H_{a,0}'(\tilde{\nu}_{ja}) \Biggr)^2 \Biggr\}\\
&&\qquad\leq K \Bigl[\max_{a-k_a+1\leq j\leq a} \bigl(\tilde H_{a,M}'(\tilde
{\nu}_{ja})-\tilde
H_{a,0}'(\tilde{\nu}_{ja}) \bigr) \Bigr]\to0\qquad \mbox{as }
a\to\infty.
\end{eqnarray*}
Since $(\sigma_a^M(k_a))^2-(\sigma_a^0(k_a))^2= (\sigma
_a^M(k_a)+\sigma_a^0(k_a) )
(\sigma_a^M(k_a)-\sigma_a^0(k_a) )$, we have
\[
\sigma_a^M(k_a)-\sigma_a^0(k_a)\to0\qquad \mbox{as } a\to\infty.
\]

\subsubsection{\texorpdfstring{Proof of Lemma \protect\ref{diff.mu}}{Proof of Lemma 4.6}}

We hope to show that
\[
\sup_{-M<t<M} \biggl|\frac{\sqrt{a}(\mu_a^t(k_a)-\mu_a^0(k_a))}
{\sigma_a^0(k_a)} \biggr|=\frac{\sqrt{a}(\mu_a^M(k_a)-\mu
_a^0(k_a))} {\sigma_a^0(k_a)}\to0\qquad \mbox{as }
a\to\infty,
\]
provided that $k_a\to\infty$, as $a\to\infty$. From the fact that
$G_{a,M}^{-1}(1-e^{-\tilde{\nu}_{ia}})
-G_{a,0}^{-1}(1-e^{-\tilde{\nu}_{ia}})$ is increasing in $i$ and
Taylor's expansion, we have
\begin{eqnarray*}
\mu_a^M(k_a)-\mu_a^0(k_a)&=&\frac{1}{a}\sum_{i=a-k_a+1}^a
\bigl(G_{a,M}^{-1}(1-e^{-\tilde{\nu}_{ia}})
-G_{a,0}^{-1}(1-e^{-\tilde{\nu}_{ia}}) \bigr)\\
&\leq&\frac{k_a}{a} \bigl(G_{a,M}^{-1}(1-e^{-\tilde{\nu}_{aa}})
-G_{a,0}^{-1}(1-e^{-\tilde{\nu}_{aa}}) \bigr)\\
&=&\frac{k_a}{a}\cdot
O \biggl(\frac{M^2}{a}G_{a,M}^{-1}(1-e^{-\tilde{\nu}_{aa}}) \biggr).
\end{eqnarray*}
Note that the last equality is justified by the similar argument of the
proof of Lemma \ref{rel.chi-nonchi}(3).
Applying the same argument of the proof of Lemma \ref
{rel.chi-nonchi}(1), it follows that
\[
G_{a,M}^{-1}(1-e^{-\tilde{\nu}_{aa}})\leq\bigl(2\log(a+1)-2\log
\bigl(\sqrt{\pi/2} \bigr)+2\log\bigl(
e^{M/(2\sqrt{a})}+e^{-M/(2\sqrt{a})} \bigr) \bigr)^2.
\]
Suppose first that $k_a/a\to0$ as $a\to\infty$. Since $\sqrt
{a/k_a}\sigma_a^0(k_a)>0$, it is enough to show
that
%
%
\begin{equation}\label{lem4.6.1-1}
\frac{a}{\sqrt{k_a}}\bigl(\mu_a^M(k_a)-\mu_a^0(k_a)\bigr)\to0\qquad
\mbox{as } a\to\infty.
\end{equation}
Since $\sqrt{\frac{k_a}{a}}\cdot\frac{G_{a,M}^{-1}(1-e^{-\tilde
{\nu}_{aa}})}{\sqrt{a}}\to0$ as $a\to\infty$,
(\ref{lem4.6.1-1}) is satisfied. Next, we suppose that for some
$0<r\leq1$, $k_a/a\to r$ as $a\to\infty$.
Then, $\sigma_a^0(k_a)>0$, so we need to prove that
%
%
\begin{equation}\label{lem4.6.1-2}
\sqrt{a}\bigl(\mu_a^M(k_a)-\mu_a^0(k_a)\bigr)\to0\qquad \mbox{as } a\to
\infty.
\end{equation}
Since $\frac{k_a}{a}\cdot\frac{G_{a,M}^{-1}(1-e^{-\tilde{\nu
}_{aa}})}{\sqrt{a}}\to0$ as $a\to\infty$,
(\ref{lem4.6.1-2}) is also satisfied.

\subsubsection{\texorpdfstring{Proof of Lemma
\protect\ref{second.lem}}{Proof of Lemma 4.7}}

Suppose that $k_a/a\to r$, $0\leq r\leq1$, and $k_a\to\infty$, as
$a\to\infty$. Then
\begin{eqnarray*}
\mu_a^0(k_a)&=&\frac{1}{a}\sum_{i=a-k_a+1}^a\tilde{H}_{a,0}(\tilde
{\nu}_{ia}) \simeq\frac{1}{a}\sum_{i=1}^a
I \biggl(\frac{i}{a+1}>\frac{a-k_a}{a+1} \biggr)G_{a,0}^{-1}
\biggl(\frac{i}{a+1} \biggr)\\
&\to&\int_0^1I(t>1-r)G_{a,0}^{-1}(t)\,dt =\int
_{G_{a,0}^{-1}(1-r)}^{\infty}ug_{a,0}(u)\,du\qquad
\mbox{as } a\to\infty.
\end{eqnarray*}
Note that if $r=1$, $\mu_a^0(k_a)\to1$, as $a\to\infty$. Also, we have
\begin{eqnarray*}
(\sigma_a^0(k_a) )^2&\simeq&\frac{1}{a^2}\sum
_{j=1}^a\sum_{l=1}^a \biggl\{
I \biggl(\frac{j}{a+1}>\frac{a-k_a}{a+1} \biggr)I \biggl(\frac
{l}{a+1}>\frac{a-k_a}{a+1} \biggr)\\
&&\hspace*{49.7pt}{}\times \biggl(1-\frac{j}{a+1} \biggr) \biggl(1-\frac{l}{a+1}
\biggr)\\
&&\hspace*{49.7pt}{}\times \frac{\min\{j/(a+1),l/(a+1) \}}
{1-\min\{j/(a+1),l/(a+1) \}}\\
&&\hspace*{49.7pt}{}\times \frac{1}{g_{a,0} (G_{a,0}^{-1} (j/(a+1) ) )}
\frac{1}{g_{a,0} (G_{a,0}^{-1} (l/(a+1) ) )}
\biggr\}\\
&\to&\int_0^1\int_0^1I(t>1-r)I(s>1-r)\bigl(\min(t,s)-ts\bigr)\\
&&\hspace*{28.4pt}{}\times \frac
{1}{g_{a,0} (G_{a,0}^{-1}(t) )}
\frac{1}{g_{a,0} (G_{a,0}^{-1}(s) )}\,dt\,
ds\\
&=&\int_0^1\int_0^1I(t>1-r)I(s>1-r)\bigl(\min(t,s)-ts\bigr)\,
dG_{a,0}^{-1}(t)\,dG_{a,0}^{-1}(s).
\end{eqnarray*}
Note that if $r=1$, $\sigma_a^0(k_a)\to\sqrt{2}$, as $a\to\infty$.

\subsection{\texorpdfstring{Proof of Theorem
\protect\ref{anova.order.thres}}{Proof of Theorem 4.3}}

From Theorem \ref{anova.t.random}, Lemmas \ref{ratio.std}, \ref
{diff.mu}, \ref{second.lem} and Slutsky's
theorem, it follows that
\begin{eqnarray*}
\widetilde T_L(k_a)&=&\frac{\widehat
T_L(k_a)-a\mu_a^0(k_a)}{\sqrt{a}\sigma_a^0(k_a)}\\
&=&s^2\frac{\widehat{\sigma}_a(k_a)}{\sigma_a^0(k_a)}\widehat
T_L^*(k_a)+s^2\frac{\sqrt{a} (\widehat{\mu}_a(k_a)-\mu
_a^0(k_a) )}{\sigma_a^0(k_a)}
+\frac{\mu_a^0(k_a)}{\sigma_a^0(k_a)}\sqrt{a} (s^2-1 )\\
&\stackrel{d}{\to}& N \biggl(0,1+\frac{2\mu_r^2}{\sigma
_r^2(n-1)} \biggr)\qquad \mbox{as } a\to\infty.
\end{eqnarray*}
\end{appendix}

\section*{Acknowledgments}
We are grateful to the Associate Editor and a referee for many helpful
comments that led to substantial improvement of the manuscript.

\begin{supplement}
\stitle{Supplement to ``Order Thresholding''}
\slink[doi]{10.1214/09-AOS782SUPP}
\sdatatype{.pdf}
\sdescription{We prove Theorems \ref{sing.seq.3.1}, \ref{anova.thm.4.1}, \ref{anova.t.random}
and \ref{anova.order.thres} of the paper
``Order Thresholding.''  A number of auxiliary results that are needed
for these proofs are also stated and proved.}
\end{supplement}

\printaddresses


\begin{thebibliography}{99}

\bibitem[\protect\citeauthoryear{Akritas and Papadatos}{2004}]{AP04}
\textsc{Akritas, M. G.} and \textsc{Papadatos, N.} (2004).
Heteroscedastic one-way ANOVA and lack-of-fit tests.
\textit{J. Amer. Statist. Assoc.} \textbf{99}
368--382.
\MR{2062823}

\bibitem[\protect\citeauthoryear{Beran}{2004}]{B04}
\textsc{Beran, R.} (2004). Hybrid shrinkage estimators using penalty
bases for the ordinal one-way layout.
\textit{Ann. Statist.} \textbf{32} 2532--2558.
\MR{2153994}

\bibitem[\protect\citeauthoryear{Bickel}{1967}]{B67}
\textsc{Bickel, P. J.} (1967). Some contributions to the theory of
order statistics. In \textit{Proc. Fifth Berkeley Symp. Math. Statist.
Probab.} \textbf{1} 575--591.
Univ. California Press, Berkeley, CA.
\MR{0216701}

\bibitem[\protect\citeauthoryear{Chernoff, Gastwirth and Johns}{1967}]{CGJ67}
\textsc{Chernoff, H.}, \textsc{Gastwirth, J. L.} and \textsc{Johns,
M. V.} (1967). Asymptotic distribution of
linear combinations of functions of order statistics with applications
to estimation. \textit{Ann. Math. Statist.} \textbf{38} 52--72.
\MR{0203874}

\bibitem[\protect\citeauthoryear{David and Nagaraja}{2003}]{DN03}
\textsc{David, H. A.} and \textsc{Nagaraja, H. N.} (2003). \textit
{Order Statistics}. Wiley, New York.
\MR{1994955}

\bibitem[\protect\citeauthoryear{Donoho and Johnstone}{1994}]{DJ94}
\textsc{Donoho, D. L.} and \textsc{Johnstone, I. M.} (1994). Ideal
spatial adaptation by wavelet shrinkage.
\textit{Biometrika} \textbf{81} 425--455.
\MR{1311089}

\bibitem[\protect\citeauthoryear{Efron et al.}{2001}]{Efronetal01}
\textsc{Efron, B.}, \textsc{Tibshirani, R.}, \textsc{Storey, J. D.}
and \textsc{Tusher, V.} (2001). Empirical
Bayes analysis of a microarray experiment.
\textit{J. Amer. Statist. Assoc.}
\textbf{96} 1151--1160.
\MR{1946571}

\bibitem[\protect\citeauthoryear{Fan}{1996}]{F96}
\textsc{Fan, J.} (1996). Test of significance based on wavelet
thresholding and Neyman's truncation.
\textit{J. Amer. Statist. Assoc.} \textbf{91}
674--688.
\MR{1395735}

\bibitem[\protect\citeauthoryear{Fan and Lin}{1998}]{FL98}
\textsc{Fan, J.} and \textsc{Lin, S. K.} (1998). Test of significance
when data are curves. \textit{J. Amer. Statist. Assoc.}
\textbf{93} 1007--1021.
\MR{1649196}

\bibitem[\protect\citeauthoryear{Galambos}{1995}]{G95}
\textsc{Galambos, J.} (1995). \textit{Advanced Probability Theory}.
Dekker, New York.
\MR{1350792}

\bibitem[\protect\citeauthoryear{Hall}{1978}]{H78}
\textsc{Hall, P.} (1978). Representations and limit theorems for
extreme value distributions. \textit{J.
Appl. Probab.} \textbf{15} 639--644.
\MR{0494433}

\bibitem[\protect\citeauthoryear{Johnstone and Silverman}{2004}]{JS04}
\textsc{Johnstone, I. M.} and \textsc{Silverman, B. W.} (2004).
Needles and straw in haystacks: Empirical Bayes estimates of
possibly sparse sequences. \textit{Ann. Statist.}
\textbf{32} 1594--1649.
\MR{2089135}


\bibitem[\protect\citeauthoryear{Kim and Akritas}{2010}]{KA09}
\textsc{Kim, M. H.} and \textsc{Akritas, M. G.} (2010). Supplement to ``Order
thresholding.''
DOI: \href{http://dx.doi.org/10.1214/09-AOS782SUPP}{10.1214/09-AOS782SUPP}.

\bibitem[\protect\citeauthoryear{Nagaraja}{1982}]{N82}
\textsc{Nagaraja, H. N.} (1982). Some nondegenerate limit laws for the
selection differential. \textit{Ann. Statist.} \textbf{10} 1306--1310.
\MR{0673667}

\bibitem[\protect\citeauthoryear{Neyman}{1937}]{N37}
\textsc{Neyman, J.} (1937). Smooth test for goodness of fit.
\textit{Skandinavisk Aktuarietidskrift}
\textbf{20} 149--199.

\bibitem[\protect\citeauthoryear{Shorack}{1969}]{Shorack69}
\textsc{Shorack, G. R.} (1969). Asymptotic normality of linear
combinations of functions of order statistics.
\textit{Ann. Math. Statist.} \textbf{40} 2041--2050.
\MR{0253457}


\bibitem[\protect\citeauthoryear{Simes}{1986}]{S86}
\textsc{Simes, R. J.} (1986). An improved Bonferroni procedure for
multiple tests of significance.
\textit{Biometrika} \textbf{73} 751--754.
\MR{0897872}

\bibitem[\protect\citeauthoryear{Spokoiny}{1996}]{S96}
\textsc{Spokoiny, V. G.} (1996). Adaptive hypothesis testing using
wavelets. \textit{Ann. Statist.}
\textbf{24} 2477--2498.
\MR{1425962}

\bibitem[\protect\citeauthoryear{Stigler}{1969}]{Stigler69}
\textsc{Stigler, S. M.} (1969). Linear functions of order statistics.
\textit{Ann. Math. Statist.} \textbf{40} 770--788.
\MR{0264822}

\bibitem[\protect\citeauthoryear{Stigler}{1973}]{S73}
\textsc{Stigler, S. M.} (1973). The asymptotic distribution of the
trimmed mean. \textit{Ann. Statist.} \textbf{1} 472--477.
\MR{0359134}

\bibitem[\protect\citeauthoryear{Storey}{2002}]{S02}
\textsc{Storey, J. D.} (2002). A direct approach to false discovery
rates. \textit{J. R. Stat. Soc. Ser. B Stat. Methodol.} \textbf{64} 479--498.
\MR{1924302}

\bibitem[\protect\citeauthoryear{Storey}{2003}]{S03}
\textsc{Storey, J. D.} (2003). The positive false discovery rate:
A Bayesian interpretation and the $q$-value.
\textit{Ann. Statist.} \textbf{31} 2013--2035.
\MR{2036398}

\end{thebibliography}
\end{document}